\input ./style/arxiv-general.cfg
\documentclass[aos,MSNbibl,seceqn,dvips]{arximspdf}
\makeatletter
   \@ifpackageloaded{graphicx}{}{\usepackage{graphicx}}
\makeatother
\doi{10.1214/15-AOS1354}% Updated by VTEXPTS2LaTeX.exe, 09.07.2015
\volume{43}
\issue{6}
\pubyear{2015}
\firstpage{2624}
\lastpage{2652}
\docsubty{FLA}

\makeatletter
\def\afrac#1#2{#1/(#2)}

\newcommand{\lleft}{\left}
\newcommand{\rrvert}{\vert}
\newcommand{\rright}{\right}
\newcommand{\rrVert}{\Vert}
\newcommand{\llvert}{\vert}
\newcommand{\llVert}{\Vert}
\renewcommand{\mid}{|}
\newcommand{\argmin}{\mathop{\operatorname{argmin}}}
\newtheorem{lemma}{Lemma}[section]
\newtheorem{proposition}{Proposition}[section]
\newtheorem{teo}{Theorem}[section]
\makeatother

\begin{document}
\begin{frontmatter}

\title{Rate-optimal graphon estimation}
\runtitle{Rate-optimal graphon estimation}

\begin{aug}
% Corresponding author: Harrison Zhou - huibin.zhou@yale.edu% Updated by VTEXPTS2LaTeX.exe, 09.07.2015 16:25
%by VTEXPTS2LaTeX.exe, 09.07.2015 15:21
\author[A]{\fnms{Chao}~\snm{Gao}\ead[label=e1]{chao.gao@yale.edu}},
\author[A]{\fnms{Yu}~\snm{Lu}\ead[label=e2]{yu.lu@yale.edu}}
\and
\author[A]{\fnms{Harrison H.}~\snm{Zhou}\corref{}\ead[label=e3]{huibin.zhou@yale.edu}\ead[label=u1,url]{http://www.stat.yale.edu/\textasciitilde hz68/}}
\runauthor{C. Gao, Y. Lu and H. H. Zhou}
\affiliation{Yale University}
%\dedicated{}
\address[A]{Department of Statistics\\
Yale University\\
New Haven, Connecticut 06511\\
USA\\
\printead{e1}\\
\phantom{E-mail: }\printead*{e2}\\
\phantom{E-mail: }\printead*{e3}\\
\printead{u1}}
\end{aug}

% HISTORY:
%
\received{\smonth{10} \syear{2014}}% Updated by VTEXPTS2LaTeX.exe,
%09.07.2015 15:21
%
\revised{\smonth{6} \syear{2015}}% Updated by VTEXPTS2LaTeX.exe,
%09.07.2015 15:21

% ABSTRACT
%
\begin{abstract}
Network analysis is becoming one of the most active research areas
in statistics. Significant advances have been made recently on
developing theories, methodologies and algorithms for analyzing
networks. However, there has been little fundamental study on optimal
estimation. In this paper, we establish optimal rate of convergence for
graphon estimation. For the stochastic block model with $k$ clusters,
we show that the optimal rate under the mean squared error is
$n^{-1}\log k+k^2/n^2$. The minimax upper bound improves the existing
results in literature through a technique of solving a quadratic
equation. When $k\leq\sqrt{n\log n}$, as the number of the cluster
$k$ grows, the minimax rate grows slowly with only a logarithmic order
$n^{-1}\log k$. A key step to establish the lower bound is to construct
a novel subset of the parameter space and then apply Fano's lemma, from
which we see a clear distinction of the nonparametric graphon
estimation problem from classical nonparametric regression, due to the
lack of identifiability of the order of nodes in exchangeable random
graph models. As an immediate application, we consider nonparametric
graphon estimation in a H\"{o}lder class with smoothness $\alpha$.
When the smoothness $\alpha\geq1$, the optimal rate of convergence is
$n^{-1}\log n$, independent of $\alpha$, while for $\alpha\in(0,1)$,
the rate is $n^{-2\alpha/(\alpha+1)}$, which is, to our surprise,
identical to the classical nonparametric rate.
\end{abstract}

% KEYWORDS
\end{frontmatter}

%s1 #&#
\section{Introduction}

Network analysis \cite{goldenberg10} has gained considerable research
interests in both theories \cite{bickel09} and applications \cite
{wasserman94,girvan02}. A lot of recent work has been focusing on
studying networks from a nonparametric perspective \cite{bickel09},
following the deep advancement in exchangeable arrays \cite
{aldous81,hoover79,kallenberg89,diaconis07}. In this paper, we study
the fundamental limits in estimating the underlying generating
mechanism of network models, called graphon. Though various algorithms
have been proposed and analyzed \cite
{chatterjee12,olhede13,wolfe13,airoldi13,chan14}, it is not clear
whether the convergence rates obtained in these works can be improved,
and not clear what the differences and connections are between
nonparametric graphon estimation and classical nonparametric regression.
The results obtained in this paper provide answers to those questions.
We found many existing results in literature are not sharp.
Nonparametric graphon estimation can be seen as nonparametric
regression without knowing design. When the smoothness of the graphon
is small, the minimax rate of graphon estimation is identical to that
of nonparametric regression. This is surprising, since graphon
estimation seems to be a more difficult problem, for which the design
is not observed. When the smoothness is high, we show that the minimax
rate does not depend on the smoothness anymore, which provides a clear
distinction between nonparametric graphon estimation and nonparametric
regression.

We consider an undirected graph of $n$ nodes. The connectivity can be
encoded by an adjacency matrix $\{A_{ij}\}$ taking values in $\{0,1\}
^{n\times n}$. The value of $A_{ij}$ stands for the presence or the
absence of an edge between the $i$th and the $j$th nodes. The model in
this paper is
$A_{ij}=A_{ji}\sim\operatorname{Bernoulli}(\theta_{ij})$ for $1\leq
j< i\leq n$,
where
%
%e1.1 #&#
%
\begin{equation}
\theta_{ij}=f(\xi_i,\xi_j),\qquad i\neq j
\in[n]. \label{eqgraphon}
\end{equation}
The sequence $\{\xi_i\}$ are random variables sampled from a
distribution $\mathbb{P}_{\xi}$ supported on $[0,1]^n$. A common
choice for the probability $\mathbb{P}_{\xi}$ is i.i.d. uniform
distribution on $[0,1]$. In this paper, we allow $\mathbb{P}_{\xi}$
to be any distribution, so that the model (\ref{eqgraphon}) is
studied to its full generality. Given $\{\xi_i\}$, we assume $\{
A_{ij}\}$ are independent for $1\leq j< i\leq n$, and adopt the
convention that $A_{ii}=0$ for each $i\in[n]$. The nonparametric model
(\ref{eqgraphon}) is inspired by the advancement of graph limit
theory \cite{lovasz06,diaconis07,lovasz12}. The function $f(x,y)$,
which is assumed to be symmetric, is called graphon. This concept plays
a significant role in network analysis. Since graphon is an object
independent of the network size $n$, it gives a natural criterion to
compare networks of different sizes. Moreover, model based prediction
and testing can be done through graphon \cite{lloyd13}. Besides
nonparametric models, various parametric models
have been proposed on the matrix $\{\theta_{ij}\}$ to capture
different aspects of the network \cite
{holland81,holland83,nowicki01,newman07,handcock07,hoff08,airoldi09,karrer11}.

The model (\ref{eqgraphon}) has a close relation to the classical
nonparametric regression problem. We may view the setting (\ref
{eqgraphon}) as modeling the mean of $A_{ij}$ by a regression function
$f(\xi_i,\xi_j)$ with design $\{(\xi_i,\xi_j)\}$. In a regression
problem, the design points $\{(\xi_i,\xi_j)\}$ are observed, and the
function $f$ is estimated from the pair $\{(\xi_i,\xi_j), A_{ij}\}$.
In contrast, in the graphon estimation setting, $\{(\xi_i,\xi_j)\}$
are latent random variables, and $f$ can only be estimated from the
response $\{A_{ij}\}$. This causes an identifiability problem, because
without observing the design, there is no way to associate the value of
$f(x,y)$ with $(x,y)$. In this paper, we consider the following loss function:
\[
\frac{1}{n^2}\sum_{i,j\in[n]}(\hat{
\theta}_{ij}-\theta_{ij})^2
\]
to overcome the identifiability issue.
This is identical to the loss function widely used in the classical
nonparametric regression problem with the form
\[
\frac{1}{n^2}\sum_{i,j\in[n]} \bigl(\hat{f}(
\xi_i,\xi_j)-f(\xi _i,\xi_j)
\bigr)^2. %
\]
Even without observing the design $\{(\xi_i,\xi_j)\}$, it is still
possible to estimate the matrix $\{\theta_{ij}\}$ by exploiting its
underlying structure modeled by (\ref{eqgraphon}).

We first consider $\{\theta_{ij}\}$ of a block structure. This
stochastic block model, proposed by \cite{holland83}, is serving as a
standard data generating process in network community detection problem
\cite{bickel09,rohe11,amini13,joseph13,lei13,cai14}. We denote the
parameter space for $\{\theta_{ij}\}$ by $\Theta_k$, where $k$ is the
number of clusters in the stochastic block model. In total, there are
an order of $k^2$ number of blocks in $\{\theta_{ij}\}$. The value of
$\theta_{ij}$ only depends on the clusters that the $i$th and the
$j$th nodes belong to. The exact definition of $\Theta_k$ is given in
Section~\ref{secSBM}. For this setting, the minimax rate for
estimating the matrix $\{\theta_{ij}\}$ is as follows.

%th1.1 #&#
%
\begin{teo} \label{teominimax2}
Under the stochastic block model, we have
\[
\inf_{\hat{\theta}}\sup_{\theta\in\Theta_k}\mathbb{E} \biggl\{
\frac{1}{n^2}\sum_{i,j\in[n]}(\hat{\theta}_{ij}-
\theta _{ij})^2 \biggr\}\asymp\frac{k^2}{n^2}+
\frac{\log k}{n},
\]
for any $1\leq k\leq n$.
\end{teo}

The convergence rate has two terms. The first term $k^2/n^2$ is due to
the fact that we need to estimate an order of $k^2$ number of unknown
parameters with an order of $n^2$ number of observations. The second
term $n^{-1}\log k$, which we coin as the clustering rate, is the error
induced by the lack of identifiability of the order of nodes in
exchangeable random graph models. Namely, it
is resulted from the unknown clustering structure of the $n$ nodes.
This term grows logarithmically as the number of clusters $k$
increases, which is different from what is obtained in literature \cite
{chatterjee12} based on lower rank matrix estimation.

We also study the minimax rate of estimating $\{\theta_{ij}\}$ modeled
by the relation (\ref{eqgraphon}) with $f$ belonging to a H\"{o}lder
class $\mathcal{F}_{\alpha}(M)$ with smoothness $\alpha$. The class
$\mathcal{F}_{\alpha}(M)$ is rigorously defined in Section~\ref
{secgraphon}. The result is stated in the following theorem.

%th1.2 #&#
%
\begin{teo} \label{teominimax1}
Consider the H\"{o}lder class $\mathcal{F}_{\alpha}(M)$, defined in
Section~\ref{secgraphon}. We have
\[
\inf_{\hat{\theta}}\sup_{f\in\mathcal{F}_{\alpha}(M)}\sup
_{\xi
\sim\mathbb{P}_{\xi}}\mathbb{E} \biggl\{\frac{1}{n^2}\sum
_{i,j\in [n]}(\hat{\theta}_{ij}-\theta_{ij})^2
\biggr\}\asymp\cases{ n^{-\afrac{2\alpha}{\alpha+1}}, &\quad$0<\alpha<1$,
\vspace*{5pt}\cr
\displaystyle
\frac{\log n}{n}, &\quad$\alpha\geq1$,} %
\]
where the expectation is jointly over $\{A_{ij}\}$ and $\{\xi_i\}$.
\end{teo}

The approximation of piecewise block function to an $\alpha$-smooth
graphon $f$ yields an additional error at the order of $k^{-2\alpha}$
(see Lemma~\ref{lembias}). In view of the minimax rate in Theorem
\ref{teominimax2}, picking the best $k$ to trade off the sum of the
three terms
$k^{-2\alpha}$, $k^2/n^2$, and $n^{-1}\log k$
gives the minimax rate in Theorem~\ref{teominimax1}.

The minimax rate reveals a new phenomenon in nonparametric estimation.
When the smoothness parameter $\alpha$ is smaller than $1$, the
optimal rate of convergence is the typical\vspace*{1pt} nonparametric rate. Note
that the typical nonparametric rate is $N^{-\afrac{2\alpha}{2\alpha
+d}}$ \cite{tsybakov09}, where $N$ is the number of observations and
$d$ is the function dimension. Here, we are in a two-dimensional
setting with number of observations $N\asymp n^2$ and dimension $d=2$.
Then the corresponding rate is $N^{-\afrac{2\alpha}{2\alpha+d}}\asymp
n^{-\afrac{2\alpha}{\alpha+1}}$. Surprisingly, in Theorem~\ref
{teominimax1} for the regime $\alpha\in(0,1)$, we get the exact same
nonparametric minimax rate, though we are not given the knowledge of
the design $\{(\xi_i,\xi_j)\}$. The cost of not observing the design
is reflected in the case with $\alpha\geq1$. In this regime, the
smoothness of the function does not help improve the rate anymore. The
minimax rate is dominated by $n^{-1}\log n$, which is essentially
contributed by the logarithmic cardinality of the set of all possible
assignments of $n$ nodes to $k$ clusters. A distinguished feature of
Theorem~\ref{teominimax1} to note is that we do not impose any
assumption on the distribution $\mathbb{P}_{\xi}$.

To prove Theorems~\ref{teominimax2}~and~\ref{teominimax1},
we develop a novel lower bound argument (see Sections~\ref{secflower}~and~\ref{secpflower}), which allows us to correctly obtain
the packing number of all possible assignments. The packing number
characterizes the difficulty brought by the ignorance of the design $\{
(\xi_i,\xi_j)\}$ in the graphon model or the ignorance of clustering
structure in the stochastic block model. Such argument may be of
independent interest, and we expect its future applications in deriving
minimax rates of other network estimation problems.

Our work on optimal graphon estimation is closely connected to a
growing literature on nonparametric network analysis. For estimating
the matrix $\{\theta_{ij}\}$ of stochastic block model, \cite
{chatterjee12} viewed $\{\theta_{ij}\}$ as a rank-$k$ matrix and
applied singular value thresholding on the adjacency matrix. The
convergence rate obtained is $\sqrt{k/n}$, which is not optimal
compared with the rate $n^{-1}\log k+k^2/n^2$ in Theorem~\ref
{teominimax2}. For nonparametric graphon estimation,
\cite{wolfe13} considered estimating $f$ in a H\"{o}lder class with
smoothness $\alpha$ and obtained the rate $\sqrt{n^{-\alpha/2} \log
n}$ under a closely related loss function. The work by \cite{chan14}
obtained the rate $n^{-1}\log n$ for estimating a Lipschitz $f$, but
they imposed strong assumptions on $f$. Namely, they assumed
$L_2\llvert  x-y\rrvert \leq\llvert  g(x)-g(y)\rrvert
\leq L_1\llvert  x-y\rrvert $ for some constants $L_1,L_2$,
with $g(x)=\int_0^1f(x,y)\,dy$. Note that this condition
excludes the stochastic block model, for which $g(x)-g(y)=0$ when
different $x$ and $y$ are in the same cluster. Local asymptotic
normality for stochastic block model was established in \cite
{bickel13}. A~method of moment via tensor decomposition was proposed by
\cite{anandkumar13}.\looseness=-1

\subsection*{Organization} The paper is organized as follows. In
Section~\ref{secmain}, we state the main results of the paper,
including both upper and lower bounds for stochastic block model and
nonparametric graphon estimation. Section~\ref{secdisc} is a
discussion section, where we discuss possible generalization of the
model, relation to nonparametric regression without knowing design and
lower bound techniques used in network analysis. The main body of the
technical proofs are presented in Section~\ref{secproof}, and the
remaining proofs are stated in the supplementary material \cite{supp}.

\subsection*{Notation} For any positive integer $d$, we use
$[d]$ to denote the set $\{1,2,\ldots,\break d\}$. For any $a,b\in\mathbb
{R}$, let $a\vee b=\max(a,b)$ and $a\wedge b=\min(a,b)$. The floor
function ${ \lfloor{a}  \rfloor}$ is the largest integer
no greater than $a$, and
the ceiling function $\lceil a \rceil$ is the smallest integer no less than
$a$. For any two positive sequences $\{a_n\}$ and $\{b_n\}$, $a_n\asymp
b_n$ means there exists a constant $C>0$ independent of $n$, such that
$C^{-1}b_n\leq a_n\leq C b_n$ for all $n$. For any $\{a_{ij}\},\{
b_{ij}\}\in\mathbb{R}^{n\times n}$, we denote the $\ell_2$ norm by
$\llVert  a\rrVert =\sqrt{\sum_{i,j\in[n]}a_{ij}^2}$ and the
inner product by
$ \langle a, b  \rangle=\sum_{i,j\in[n]}a_{ij}b_{ij}$.
Given any set $S$, $\llvert  S\rrvert $ denotes its
cardinality, and $\mathbb{I}\{
x\in S\}$ stands for the indicator function which takes value $1$ when
$x\in S$ and takes value $0$ when $x\notin S$. For a metric space
$(T,\rho)$, the covering number $\mathcal{N}(\varepsilon,T,\rho)$ is
the smallest number of balls with radius $\varepsilon$ and centers in $T$
to cover $T$, and the packing number $\mathcal{M}(\varepsilon,T,\rho)$
is the largest number of points in $T$ that are at least $\varepsilon$
away from each other. The symbols $\mathbb{P}$ and $\mathbb{E}$ stand
for generic probability and expectation, whenever the distribution is
clear from the context.

%s2 #&#
\section{Main results} \label{secmain}

In this section, we present the main results of the paper. We first
introduce the estimation procedure in Section~\ref{secmethod}. The
minimax rates of stochastic block and nonparametric graphon estimation
are stated in Sections~\ref{secSBM}~and~\ref{secgraphon},
respectively.

%s2.1 #&#
\subsection{Methodology} \label{secmethod}

We are going to propose an estimator for both stochastic block model
and nonparametric graphon estimation under H\"{o}lder smoothness.
To introduce the estimator, let us define the set $\mathcal
{Z}_{n,k}= \{z: [n]\rightarrow[k] \}$ to be the collection
of all possible mappings from $[n]$ to $[k]$ with some integers $n$ and $k$.
Given a $z\in\mathcal{Z}_{n,k}$, the sets $\{z^{-1}(a):a\in[k]\}$
form a partition of $[n]$, in the sense that $\bigcup_{a\in
[k]}z^{-1}(a)=[n]$ and $z^{-1}(a)\cap z^{-1}(b)=\varnothing$ for any
$a\neq b\in[k]$. In other words, $z$ defines a clustering structure on
the $n$ nodes. It is easy to see that the cardinality of $\mathcal
{Z}_{n,k}$ is $k^n$.
Given a matrix $\{\eta_{ij}\}\in\mathbb{R}^{n\times n}$, and a
partition function $z\in\mathcal{Z}_{n,k}$, we use the following
notation to denote the block average on the set $z^{-1}(a)\times
z^{-1}(b)$. That is,
%
%e2.1 #&#
%
\begin{equation}
\bar{\eta}_{ab}(z)=\frac{1}{\llvert  z^{-1}(a)\rrvert \llvert  z^{-1}(b)\rrvert }\sum_{i\in
z^{-1}(a)}
\sum_{j\in z^{-1}(b)}\eta_{ij}\qquad\mbox{for } a\neq b
\in[k],\label{eqaveab}
\end{equation}
and when $\llvert  z^{-1}(a)\rrvert >1$,
%
%e2.2 #&#
%
\begin{equation}
\bar{\eta}_{aa}(z)=\frac{1}{\llvert  z^{-1}(a)\rrvert (\llvert  z^{-1}(a)\rrvert -1)}\sum_{i\neq
j\in z^{-1}(a)}
\eta_{ij}\qquad\mbox{for } a\in[k].\label{eqaveaa}
\end{equation}
For any $Q=\{Q_{ab}\}\in\mathbb{R}^{k\times k}$ and $z\in\mathcal
{Z}_{n,k}$, define the objective function
\[
L(Q,z)=\sum_{a,b\in[k]}\mathop{\sum
_{(i,j)\in z^{-1}(a)\times
z^{-1}(b)}}_{i\neq j}(A_{ij}-Q_{ab})^2.
\]
For any optimizer of the objective function,
%
%e2.3 #&#
%
\begin{equation}
(\hat{Q},\hat{z})\in\argmin_{Q\in\mathbb{R}^{k\times k},z\in
\mathcal{Z}_{n,k}}L(Q,z), \label{eqestimator}
\end{equation}
the estimator of $\theta_{ij}$ is defined as
%
%e2.4 #&#
%
\begin{equation}
\hat{\theta}_{ij}=\hat{Q}_{\hat{z}(i)\hat{z}(j)},\qquad i>j, \label
{eqdefthetahat}
\end{equation}
and $\hat{\theta}_{ij}=\hat{\theta}_{ji}$ for $i<j$. Set the
diagonal element by $\hat{\theta}_{ii}=0$. The procedure (\ref
{eqdefthetahat}) can be understood as first clustering the data by an
estimated $\hat{z}$ and then estimating the model parameters via block
averages. By the least squares formulation, it is easy to observe the
following property.
%pr2.1 #&#

\begin{proposition} \label{propestimation}
For any minimizer $(\hat{Q},\hat{z})$, the entries of $\hat{Q}$ has
representation
%
%e2.5 #&#
%
\begin{equation}
\hat{Q}_{ab}=\bar{A}_{ab}(\hat{z}), \label{eqdefQhat}
\end{equation}
for all $a,b\in[k]$.
\end{proposition}

The representation of the solution (\ref{eqdefQhat}) shows that the
estimator (\ref{eqdefthetahat}) is essentially doing a histogram
approximation after finding the optimal cluster assignment $\hat{z}\in
\mathcal{Z}_{n,k}$ according to the least squares criterion (\ref
{eqestimator}). In the classical nonparametric regression problem, it
is known that a simple histogram estimator cannot achieve optimal
convergence rate for $\alpha>1$ \cite{tsybakov09}. However, we are
going to show that this simple histogram estimator achieves optimal
rates of convergence under both stochastic block model and
nonparametric graphon estimation settings.

Similar estimators using the Bernoulli likelihood function have been
proposed and analyzed in the literature \cite
{bickel09,zhao12,wolfe13,olhede13}. Instead of using the likelihood
function of Bernoulli distribution, the least squares estimator (\ref
{eqestimator}) can be viewed as maximizing Gaussian likelihood. This
allows us to obtain optimal convergence rates with cleaner analysis.

%s2.2 #&#
\subsection{Stochastic block model} \label{secSBM}

In the stochastic block model setting, each node $i\in[n]$ is
associated with a label $a\in[k]$, indicating its cluster.
The edge $A_{ij}$ is a Bernoulli random variable with mean $\theta
_{ij}$. The value of $\theta_{ij}$ only depends on the clusters of the
$i$th and the $j$th nodes. We assume $\{\theta_{ij}\}$ is from the
following parameter space:
\begin{eqnarray*}
\Theta_k &=& \bigl\{\{\theta_{ij}\}\in[0,1]^{n\times n}:
\theta _{ii}=0, \theta_{ij}=Q_{ab}=Q_{ba}
\\
&& \mbox{for }(i,j)\in z^{-1}(a)\times z^{-1}(b)
\mbox{ for some }Q_{ab}\in[0,1]\mbox{ and } z\in
\mathcal{Z}_{n,k} \bigr\}.
\end{eqnarray*}
Namely, the partition function $z$ assigns cluster to each node, and
the value of $Q_{ab}$ measures the intensity of link between the $a$th
and the $b$th clusters.
The least squares estimator (\ref{eqestimator}) attains the following
convergence rate for estimating $\{\theta_{ij}\}$.

%th2.1 #&#
%
\begin{teo} \label{teomain2}
For any constant $C'>0$, there is a constant $C>0$ only depending on
$C'$, such that
\[
\frac{1}{n^2}\sum_{i,j\in[n]}(\hat{
\theta}_{ij}-\theta _{ij})^2\leq C \biggl(
\frac{k^2}{n^2}+\frac{\log k}{n} \biggr),
\]
with probability at least $1-\exp(-C'n\log k)$, uniformly over $\theta
\in\Theta_k$. Furthermore, we have
\[
\sup_{\theta\in\Theta_k}\mathbb{E} \biggl\{\frac{1}{n^2}\sum
_{i,j\in[n]}(\hat{\theta}_{ij}-\theta_{ij})^2
\biggr\}\leq C_1 \biggl(\frac{k^2}{n^2}+\frac{\log k}{n} \biggr),
\]
for all $k\in[n]$ with some universal constant $C_1>0$.
\end{teo}

Theorem~\ref{teomain2} characterizes different convergence rates for
$k$ in different regimes. Suppose $k\asymp n^{\delta}$ for some
$\delta\in[0,1]$. Then the convergence rate in Theorem~\ref
{teomain2} is
%
%e2.6 #&#
%
\begin{equation}
\frac{k^2}{n^2}+\frac{\log k}{n}\asymp\cases{ n^{-2}, & \quad$k=1$,
\vspace*{5pt}\cr
n^{-1}, &\quad$\delta=0$, $k\geq2$,
\vspace*{5pt}\cr
n^{-1}\log n, &
\quad$\delta\in(0,1/2]$,
\vspace*{5pt}\cr
n^{-2(1-\delta)}, &\quad$\delta\in(1/2,1]$.}\label{eqfenlei}
\end{equation}
The result completely characterizes the convergence rates for
stochastic block model with any possible number of clusters $k$.
Depending on whether $k$ is small, moderate or large, the convergence
rates behave differently.

The convergence rate, in terms of $k$, has two parts. The first part
$k^2/n^2$ is called the nonparametric rate. It is determined by the
number of parameters and the number of observations of the model. For
the stochastic block model with $k$ clusters, the number of parameters
is $k(k+1)/2\asymp k^2$ and\vspace*{1pt} the number of observations is
$n(n+1)/2\asymp n^2$. The second part $n^{-1}\log k$ is called the
clustering rate. Its presence is due to the unknown labels of the $n$
nodes. Our result shows the clustering rate is logarithmically
depending on the number of clusters $k$.
From (\ref{eqfenlei}), we observe that when $k$ is small, the
clustering rate dominates. When $k$ is large, the nonparametric rate dominates.

To show that the rate in Theorem~\ref{teomain2} cannot be improved,
we obtain the following minimax lower bound.

%th2.2 #&#
%
\begin{teo} \label{teoSBMlower}
There exists a universal constant $C>0$, such that
\[
\inf_{\hat{\theta}}\sup_{\theta\in\Theta_k}\mathbb{P} \biggl\{
\frac{1}{n^2}\sum_{i,j\in[n]}(\hat{\theta}_{ij}-
\theta _{ij})^2\geq C \biggl(\frac{k^2}{n^2}+
\frac{\log k}{n} \biggr) \biggr\} \geq0.8,
\]
and
\[
\inf_{\hat{\theta}}\sup_{\theta\in\Theta_k}\mathbb{E} \biggl\{
\frac{1}{n^2}\sum_{i,j\in[n]}(\hat{\theta}_{ij}-
\theta _{ij})^2 \biggr\}\geq C \biggl(\frac{k^2}{n^2}+
\frac{\log k}{n} \biggr),
\]
for any $k\in[n]$.
\end{teo}

The upper bound of Theorem~\ref{teomain2} and the lower bound of
Theorem~\ref{teoSBMlower} immediately imply the minimax rate in
Theorem~\ref{teominimax2}.

%s2.3 #&#
\subsection{Nonparametric graphon estimation} \label{secgraphon}

Let us proceed to nonparametric graphon estimation. For any $i\neq j$,
$A_{ij}$ is sampled from the following process:
\[
(\xi_1,\ldots,\xi_n)\sim\mathbb{P}_{\xi},\qquad
A_{ij}\mid(\xi _i,\xi _j)\sim
\operatorname{Bernoulli}(\theta_{ij})\qquad\mbox{where }
\theta_{ij}=f(\xi_i,\xi_j).
\]
For $i\in[n]$, $A_{ii}=\theta_{ii}=0$.
Conditioning on $(\xi_1,\ldots,\xi_n)$, $A_{ij}$ is independent
across $i,j\in[n]$.
To completely specify the model, we need to define the function class
of $f$ on $[0,1]^2$. Since $f$ is symmetric, we only need to specify
its value on $\mathcal{D}=\{(x,y)\in[0,1]^2: x\geq y\}$. Define the
derivative operator by
\[
\nabla_{jk}f(x,y)=\frac{\partial^{j+k}}{(\partial x)^j(\partial y)^k}f(x,y),
\]
and we adopt the convention $\nabla_{00}f(x,y)=f(x,y)$.
The H\"{o}lder norm is defined as
\begin{eqnarray*}
\llVert f\rrVert _{\mathcal{H}_{\alpha}} &=& \max_{j+k\leq
{ \lfloor{\alpha}  \rfloor}}\sup
_{x,y\in\mathcal{D}}\bigl\llvert \nabla_{jk}f(x,y)\bigr\rrvert
\\
&&{} +\max
_{j+k={ \lfloor{\alpha}  \rfloor}}\sup_{(x,y)\neq
(x',y')\in\mathcal{D}}\frac
{\llvert \nabla_{jk}f(x,y)-\nabla_{jk}f(x',y')\rrvert
}{(\llvert  x-x'\rrvert +\llvert  y-y'\rrvert )^{\alpha
-{ \lfloor{\alpha}  \rfloor}}}.
\end{eqnarray*}
The H\"{o}lder class is defined by
\[
\mathcal{H}_{\alpha}(M)= \bigl\{\| f\| _{\mathcal{H}_{\alpha}}\leq
M: f(x,y)=f(y,x)\mbox{ for }x\geq y \bigr\},
\]
where $\alpha>0$ is the smoothness parameter and $M>0$ is the size of
the class, which is assumed to be a constant. When $\alpha\in(0,1]$,
a function $f\in\mathcal{H}_{\alpha}(M)$ satisfies the Lipschitz condition
%
%e2.7 #&#
%
\begin{equation}
\bigl\llvert f(x,y)-f\bigl(x',y'\bigr)\bigr\rrvert
\leq M\bigl(\bigl\llvert x-x'\bigr\rrvert +\bigl\llvert
y-y'\bigr\rrvert \bigr)^{\alpha}, \label{eqLip}
\end{equation}
for any $(x,y), (x',y')\in\mathcal{D}$.
In the network model,
the graphon $f$ is assumed to live in the following class:
\[
\mathcal{F}_{\alpha}(M)= \bigl\{0\leq f\leq1: f\in\mathcal
{H}_{\alpha}(M) \bigr\}.
\]
We have mentioned that the convergence rate of graphon estimation is
essentially due to the stochastic block model approximation of $f$ in a
H\"{o}lder class.
This intuition is established by the following lemma, whose proof is
given in the supplementary material \cite{supp}.
%\nb{using $\tilde{z}$ or $z^*$?}
%le2.1 #&#

\begin{lemma} \label{lembias}
There exists $z^*\in\mathcal{Z}_{n,k}$, satisfying
\[
\frac{1}{n^2}\sum_{a,b\in[k]}\sum
_{\{i\neq j: z^*(i)=a, z^*(j)=b\}
} \bigl(\theta_{ij}-\bar{\theta}_{ab}
\bigl(z^*\bigr) \bigr)^2\leq CM^2 \biggl(\frac{1}{k^2}
\biggr)^{\alpha\wedge1},
\]
for some universal constant $C>0$.
\end{lemma}

The graph limit theory \cite{lovasz06} suggests $\mathbb{P}_{\xi}$
to be an i.i.d. uniform distribution on the interval $[0,1]$. For the
estimating procedure (\ref{eqestimator}) to work, we allow $\mathbb
{P}_{\xi}$ to be any distribution. The upper bound is attained over
all distributions $\mathbb{P}_{\xi}$ uniformly. Combining Lemma~\ref
{lembias} and Theorem~\ref{teomain2} in an appropriate manner, we
obtain the convergence rate for graphon estimation by the least squares
estimator (\ref{eqestimator}).

%th2.3 #&#
%
\begin{teo} \label{teomain1}
Choose $k=\lceil n^{\afrac{1}{\alpha\wedge1+1}} \rceil$. Then for any
$C'>0$, there exists a constant $C>0$ only depending on $C'$ and $M$,
such that
\[
\frac{1}{n^2}\sum_{i,j\in[n]}(\hat{
\theta}_{ij}-\theta _{ij})^2\leq C
\biggl(n^{-\afrac{2\alpha}{\alpha+1}}+\frac{\log
n}{n} \biggr),
\]
with probability at least $1-\exp(-C'n)$, uniformly over $f\in
\mathcal{F}_{\alpha}(M)$ and $\mathbb{P}_{\xi}$. Furthermore,
\[
\sup_{f\in\mathcal{F}_{\alpha}(M)}\sup_{\mathbb{P}_{\xi}}\mathbb {E} \biggl\{
\frac{1}{n^2}\sum_{i,j\in[n]}(\hat{\theta}_{ij}-
\theta _{ij})^2 \biggr\}\leq C_1
\biggl(n^{-\afrac{2\alpha}{\alpha+1}}+\frac
{\log n}{n} \biggr),
\]
for some other constant $C_1>0$ only depending on $M$. Both the
probability and the expectation are jointly over $\{A_{ij}\}$ and $\{
\xi_i\}$.
\end{teo}

Similar to Theorem~\ref{teomain2}, the convergence rate of Theorem
\ref{teomain1} has two parts. The nonparametric rate $n^{-\afrac{2\alpha}{\alpha+1}}$, and the clustering rate $n^{-1}\log n$. Note
that the clustering rates in both theorems are identical because
$n^{-1}\log n\asymp n^{-1}\log k$ under the choice $k=\lceil n^{\afrac{1}{\alpha\wedge1+1}} \rceil$. An interesting phenomenon to note is that
the smoothness index $\alpha$ only plays a role in the regime $\alpha
\in(0,1)$. The convergence rate is always dominated by $n^{-1}\log n$
when $\alpha\geq1$.

In order to show the rate of Theorem~\ref{teomain1} is optimal, we
need a lower bound over the class $\mathcal{F}_{\alpha}(M)$ and over
all $\mathbb{P}_{\xi}$. To be specific, we need to show
%
%e2.8 #&#
%
\begin{equation}
\qquad\inf_{\hat{\theta}}\sup_{f\in\mathcal{F}_{\alpha}(M)}\sup
_{\mathbb{P}_{\xi}}\mathbb{E} \biggl\{\frac{1}{n^2}\sum
_{i,j\in
[n]}(\hat{\theta}_{ij}-\theta_{ij})^2
\biggr\}\geq C \biggl(n^{-\afrac{2\alpha}{\alpha+1}}+\frac{\log n}{n} \biggr), \label{eqweakL}
\end{equation}
for some constant $C>0$.
In fact, the lower bound we obtained is stronger than (\ref{eqweakL})
in the sense that it holds for a subset of the space of probabilities
on $\{\xi_i\}$. The subset $\mathcal{P}$
requires the sampling points $\{\xi_i\}$ to well cover the interval
$[0,1]$ for $\{f(\xi_i,\xi_j)\}_{i,j\in[n]}$ to be good
representatives of the whole function $f$. For each $a\in[k]$, define
the interval
%
%e2.9 #&#
%
\begin{equation}
U_a=\biggl[\frac{a-1}{k},\frac{a}{k} \biggr). \label{eqUa}
\end{equation}
We define the distribution class by
\[
\mathcal{P}= \Biggl\{\mathbb{P}_{\xi}: \mathbb{P}_{\xi}
\Biggl(\frac
{\lambda_1n}{k}\leq\sum_{i=1}^n
\mathbb{I}\{\xi_i\in U_a\}\leq \frac{\lambda_2n}{k}
\mbox{ for any }a\in[k] \Biggr)>1-\exp \bigl(-n^{\delta}\bigr) \Biggr\},
\]
for some positive constants $\lambda_1,\lambda_2$ and some arbitrary
small constant $\delta\in(0,1)$. Namely, for each interval $U_a$, it
contains roughly $n/k$ observations. By applying standard concentration
inequality, it can be shown that the i.i.d. uniform distribution on $\{
\xi_i\}$ belongs to the class $\mathcal{P}$.

%th2.4 #&#
%
\begin{teo} \label{teographonlower}
There exists a constant $C>0$ only depending on $M,\alpha$, such that
\[
\inf_{\hat{\theta}}\sup_{f\in\mathcal{F}_{\alpha}(M)}\sup
_{\mathbb{P}_{\xi}\in\mathcal{P}}\mathbb{P} \biggl\{\frac
{1}{n^2}\sum
_{i,j\in[n]}(\hat{\theta}_{ij}-\theta_{ij})^2
\geq C \biggl(n^{-\afrac{2\alpha}{\alpha+1}}+\frac{\log n}{n} \biggr) \biggr\}\geq0.8,
\]
and
\[
\inf_{\hat{\theta}}\sup_{f\in\mathcal{F}_{\alpha}(M)}\sup
_{\mathbb{P}_{\xi}\in\mathcal{P}}\mathbb{E} \biggl\{\frac
{1}{n^2}\sum
_{i,j\in[n]}(\hat{\theta}_{ij}-\theta_{ij})^2
\biggr\} \geq C \biggl(n^{-\afrac{2\alpha}{\alpha+1}}+\frac{\log n}{n} \biggr),
\]
where the probability and expectation are jointly over $\{A_{ij}\}$ and
$\{\xi_i\}$.
\end{teo}

The proof of Theorem~\ref{teographonlower} is given in the
supplementary material \cite{supp}.
The minimax rate in Theorem~\ref{teominimax1} is an immediate
consequence of Theorems~\ref{teomain1} and~\ref{teographonlower}.

%s3 #&#
\section{Discussion} \label{secdisc}

%s3.1 #&#
\subsection{More general models} \label{secgeneral}

The results in this paper assume symmetry on the graphon $f$ and the
matrix $\{\theta_{ij}\}$. Such assumption is naturally made in the
context of network analysis. However, these results also hold under
more general models. We may consider a slightly more general version of
(\ref{eqgraphon}) as
\[
\theta_{ij}=f(\xi_i,\eta_j),\qquad1\leq i,j
\leq n,
\]
with $\{\xi_i\}$ and $\{\eta_j\}$ sampled from $\mathbb{P}_{\xi}$
and $\mathbb{P}_{\eta}$, respectively, and the function $f$ is not
necessarily symmetric. To be specific, let us redefine the H\"{o}lder
norm $\llVert \cdot\rrVert _{\mathcal{H}_{\alpha}}$ by
replacing $\mathcal{D}$
with $[0,1]^{2}$ in its original definition in Section~\ref
{secgraphon}. Then we consider the function class
\[
\mathcal{F}_{\alpha}'(M)=\bigl\{0\leq f\leq1: \llVert f
\rrVert _{\mathcal{H}_{\alpha
}}\leq M\bigr\}.
\]
The minimax rate for this class is stated in the following theorem
without proof.
%th3.1 #&#

\begin{teo}\label{teodisc1}
Consider the function class $\mathcal{F}_{\alpha}'(M)$ with $\alpha
>0$ and $M>0$. We have
\[
\inf_{\hat{\theta}}\sup_{f\in\mathcal{F}'_{\alpha}(M)} \mathop{\sup
_{\xi\sim\mathbb{P}_{\xi}}}_{\eta\sim\mathbb
{P}_{\eta}}\mathbb{E} \biggl\{\frac{1}{n^2}\sum
_{i,j\in[n]}(\hat {\theta}_{ij}-
\theta_{ij})^2 \biggr\}\asymp\cases{ n^{-\afrac{2\alpha}{\alpha+1}}, &
\quad$0<\alpha<1$,
\vspace*{5pt}\cr
\displaystyle\frac{\log n}{n}, &\quad$\alpha\geq1$,}
\]
where the expectation is jointly over $\{A_{ij}\}$, $\{\xi_i\}$ and $\{
\eta_j\}$.
\end{teo}

Similarly, we may generalize the stochastic block model by the
parameter space
\begin{eqnarray*}
\Theta_{kl}^{\mathrm{asym}} &=& \bigl\{\{\theta_{ij}\}
\in[0,1]^{n\times m}: \theta_{ij}=Q_{ab}\mbox{ for
}(i,j)\in z_1^{-1}(a)\times z_2^{-1}(b)
\\
&&\mbox{with some }Q_{ab}\in[0,1], z_1\in
\mathcal{Z}_{n,k}\mbox{ and }z_2\in\mathcal{Z}_{m,l}
\bigr\}.
\end{eqnarray*}
Such model naturally arises in the contexts of biclustering \cite
{hartigan72,mirkin98,cheng00,madeira04} and matrix organization \cite
{gavish10,coifman11,gavish12}, where symmetry of the model is not
assumed. Under such extension, we can show that a similar minimax rate
as in Theorem~\ref{teominimax2} as follows.
%th3.2 #&#

\begin{teo}\label{teodisc2}
Consider the parameter space $\Theta_{kl}^{\mathrm{asym}}$ and assume $\log
k\asymp\log l$. We have
\[
\inf_{\hat{\theta}}\sup_{\theta\in\Theta_{kl}^{\mathrm{asym}}}\mathbb {E} \biggl\{
\frac{1}{nm} \mathop{\sum_{i\in[n]}}_{j\in[m]}(
\hat{\theta}_{ij}-\theta _{ij})^2 \biggr\}\asymp
\frac{kl}{nm}+\frac{\log k}{m}+\frac{\log l}{n},
\]
for any $1\leq k\leq n$ and $1\leq l\leq m$.
\end{teo}

The lower bounds of Theorems~\ref{teodisc1}~and~\ref
{teodisc2} are directly implied by viewing the symmetric parameter
spaces as subsets of the asymmetric ones. For the upper bound, we
propose a modification of the least squares estimator in Section~\ref
{secmethod}. Consider the criterion function
\[
L^{\mathrm{asym}}(Q,z_1,z_2)=\sum
_{(a,b)\in[k]\times[l]}\sum_{(i,j)\in
z_1^{-1}(a)\times z_2^{-1}(b)}(A_{ij}-Q_{ab})^2.
\]
For any $(\hat{Q},\hat{z}_1,\hat{z}_2)\in\argmin_{Q\in\mathbb
{R}^{k\times l},z_1\in\mathcal{Z}_{n,k},z_2\in\mathcal
{Z}_{m,l}}L(Q,z_1,z_2)$, define the estimator of $\theta_{ij}$ by
\[
\hat{\theta}_{ij}=\hat{Q}_{\hat{z}_1(i)\hat{z}_2(j)}\qquad\mbox {for all } (i,j)
\in[n]\times[m].
\]
Using the same proofs of Theorems~\ref{teomain2} and~\ref
{teomain1}, we can obtain the upper bounds.

%s3.2 #&#
\subsection{Nonparametric regression without knowing design} \label
{secregression}

The graphon estimation problem is closely related to the classical
nonparametric regression problem. This section explores their
connections and differences to bring better understandings of both
problems. Namely, we study the problem of nonparametric regression
without observing the design.
First, let us consider the one-dimensional regression problem
\[
y_i=f(\xi_i)+z_i,\qquad i\in[n],
\]
where $\{\xi_i\}$ are sampled from some $\mathbb{P}_{\xi}$, and
$z_i$ are i.i.d. $N(0,1)$ variables. A~nonparametric function estimator
$\hat{f}$ estimates the function $f$ from the pairs $\{(\xi_i,y_i)\}
$. For\vspace*{1pt} H\"{o}lder class with smoothness $\alpha$, the minimax rate
under the loss $\frac{1}{n}\sum_{i\in[n]} (\hat{f}(\xi
_i)-f(\xi_i) )^2$ is at the order of $n^{-\afrac{2\alpha }{2\alpha+1}}$ \cite{tsybakov09}. However, when the design $\{\xi_i\}
$ is not observed, the minimax rate is at a constant order. To see this
fact, let us consider a closely related problem
\[
y_i=\theta_i+z_i,\qquad i\in[n],
\]
where we assume $\theta\in\Theta_2$. The parameter space $\Theta_2$
is defined as a subset of $[0,1]^n$ with $\{\theta_i\}$ that can only
take two possible values $q_1$ and $q_2$. It can be viewed as a
one-dimensional version of stochastic block model. We can show that
\[
\inf_{\hat{\theta}}\sup_{\theta\in\Theta_2}\mathbb{E} \biggl\{
\frac{1}{n}\sum_{i\in[n]}(\hat{
\theta}_i-\theta_i)^2 \biggr\} \asymp1.
\]
The upper bound is achieved by letting $\hat{\theta}_i=y_i$ for each
$i\in[n]$. To see the lower bound, we may fix $q_1=1/4$ and $q_2=1/2$.
Then the problem is reduced to $n$ independent two-point testing
problems between $N(1/4,1)$ and $N(1/2,1)$ for each $i\in[n]$. It is
easy to see that each testing problem contributes to an error at the
order of a constant, which gives the lower bound of a constant order.
This leads to a constant lower bound for the original regression
problem by using the embedding technique in the proof of Theorem~\ref
{teographonlower}, which shows that $\Theta_2$ is a smaller space
than a H\"{o}lder class on a subset of $[n]$. Thus, $1$ is also a lower
bound for the regression problem without knowing design.

In contrast to the one-dimensional problem, we can show that a
two-dimensional nonparametric regression without knowing design is more
informative. Consider
\[
y_{ij}=f(\xi_i,\xi_j)+z_{ij},\qquad
i,j\in[n],
\]
where $\{\xi_i\}$ are sampled from some $\mathbb{P}_{\xi}$, and
$z_{ij}$ are i.i.d. $N(0,1)$ variables. Let us consider the H\"{o}lder
class $\mathcal{H}_{\alpha}'(M)=\{f: \llVert  f\rrVert
_{\mathcal{H}_{\alpha
}}\leq M\}$ with H\"{o}lder norm $\llVert \cdot\rrVert
_{\mathcal{H}_{\alpha}}$
defined in Section~\ref{secgeneral}.
When the design $\{\xi_i\}$ is known, the minimax rate under the loss
$\frac{1}{n^2}\sum_{i,j\in[n]} (\hat{f}(\xi_i,\xi_j)-f(\xi
_i,\xi_j) )^2$ is at the order of $n^{-\afrac{2\alpha}{\alpha
+1}}$. When the design is unknown, the minimax rate is stated in the
following theorem.
%th3.3 #&#

\begin{teo}
Consider the H\"{o}lder class $\mathcal{H}'_{\alpha}(M)$ for $\alpha
>0$ and $M>0$. We have
\begin{eqnarray*}
&& \inf_{\hat{f}}\sup_{f\in\mathcal{H}'_{\alpha}(M)}\sup
_{\mathbb
{P}_{\xi}}\mathbb{E} \biggl\{\frac{1}{n^2}\sum
_{i,j\in[n]} \bigl(\hat{f}(\xi_i,\xi_j)-f(
\xi_i,\xi_j) \bigr)^2 \biggr\}
\\
&&\qquad \asymp \cases{
n^{-\afrac{2\alpha}{\alpha+1}}, &\quad$0<\alpha<1$,
\vspace*{5pt}\cr
\displaystyle\frac{\log n}{n}, &\quad$
\alpha\geq1$,} %
\end{eqnarray*}
where the expectation is jointly over $\{A_{ij}\}$ and $\{\xi_i\}$.
\end{teo}

The minimax rate is identical to that of Theorem~\ref{teominimax1},
which demonstrates the close relation between nonparametric graphon
estimation and nonparametric regression without knowing design.
The proof of this result is similar to the proofs of Theorems~\ref
{teomain1}~and~\ref{teographonlower}, and is omitted in the
paper. One simply needs to replace the Bernoulli analysis by the
corresponding Gaussian analysis in the proof. Compared with the rate
for one-dimensional regression without knowing design, the
two-dimensional minimax rate is more interesting. It shows that the
ignorance of design only matters when $\alpha\geq1$. For $\alpha\in
(0,1)$, the rate is exactly the same as the case when the design is known.

The main reason for the difference between the one-dimensional and the
two-dimensional problems is that the form of $\{(\xi_i,\xi_j)\}$
implicitly imposes more structure. To illustrate this point, let us
consider the following two-dimensional problem
\[
y_{ij}=f(\xi_{ij})+z_{ij},\qquad i,j\in[n],
\]
where $\xi_{ij}\in[0,1]^2$ and $\{\xi_{ij}\}$ are sampled from some
distribution. It is easy to see that this is equivalent to the
one-dimensional problem with $n^2$ observations and the minimax rate is
at the order of a constant. The form $\{(\xi_i,\xi_j)\}$ implies that
the lack of identifiability caused by the ignorance of design is only
resulted from row permutation and column permutation, and thus it is
more informative than the design $\{\xi_{ij}\}$.

%s3.3 #&#
\subsection{Lower bound for finite $k$} \label{secflower}

A key contribution of the paper lies in the proof of Theorem~\ref
{teoSBMlower}, where we establish the lower bound $k^2/n^2+n^{-1}\log
k$ (especially the $n^{-1}\log k$ part) via a novel construction. To
better understand the main idea behind the construction, we present the
analysis for a finite $k$ in this section. When $2\leq k\leq O(1)$, the
minimax rate becomes $n^{-1}$. To prove this lower bound, it is
sufficient to consider the parameter space $\Theta_k$ with $k=2$. Let
us define
\[
Q= \left[\matrix{ \displaystyle\frac{1}{2} & \displaystyle
\frac{1}{2}+\frac{c}{\sqrt
{n}}
\vspace*{5pt}\cr
\displaystyle\frac{1}{2} +
\frac{c}{\sqrt{n}} & \displaystyle\frac
{1}{2}} \right],
\]
for some $c>0$ to be determined later. Define the subspace
\[
T= \bigl\{\{\theta_{ij}\}\in[0,1]^{n\times n}: \theta
_{ij}=Q_{z(i)z(j)}\mbox{ for some }z\in\mathcal{Z}_{n,2}
\bigr\}.
\]
It is easy to see that $T\subset\Theta_2$. With a fixed $Q$, the set
$T$ has a one-to-one correspondence with $\mathcal{Z}_{n,2}$. Let us
define the collection of subsets $\mathcal{S}=\{S: S\subset[n]\}$.
For any $z\in\mathcal{Z}_{n,2}$, it induces a partition $\{
z^{-1}(1),z^{-1}(2)\}$ on the set $[n]$. This corresponds to $\{S,S^c\}
$ for some $S\in\mathcal{S}$. With this observation, we may rewrite
$T$ as
\begin{eqnarray*}
T &=& \biggl\{\{\theta_{ij}\}\in[0,1]^{n\times n}:
\theta_{ij}=\frac
{1}{2}\mbox{ for }(i,j)\in(S\times S)\cup
\bigl(S^c\times S^c\bigr),
\\
&&{} \theta_{ij} =\frac{1}{2}+\frac{c}{\sqrt{n}}\mbox{ for
}(i,j)\in\bigl(S\times S^c\bigr)\cup\bigl(S^c\times S
\bigr),\mbox{ with some }S\in \mathcal{S} \biggr\}.
\end{eqnarray*}
The subspace $T$ characterizes the difficulty of the problem due to the
ignorance of the clustering structure $\{S,S^c\}$ of the $n$ nodes.
Such difficulty is central in the estimation problem of network analysis.
We are going to use Fano's lemma (Proposition~\ref{propfano}) to
lower bound the risk. Then it is sufficient to upper bound the KL
diameter $\sup_{\theta,\theta'\in T}D(\mathbb{P}_{\theta}\|
\mathbb
{P}_{\theta'})$ and lower bound the packing number $\mathcal
{M}(\varepsilon,T,\rho)$ for some appropriate $\varepsilon$ and the metric
$\rho(\theta,\theta')=n^{-1}\llVert \theta-\theta'\rrVert
$. Using Proposition~\ref{propprobdistance}, we have
\[
\sup_{\theta,\theta'\in T}D(\mathbb{P}_{\theta}\| \mathbb
{P}_{\theta'})\leq\sup_{\theta,\theta'\in T}8\| \theta-
\theta '\| ^2\leq8c^2n.
\]
To obtain a lower bound for $\mathcal{M}(\varepsilon,T,\rho)$, note
that for $\theta,\theta'\in T$ associated with $S,S'\in\mathcal
{S}$, we have
\[
n^2\rho^2\bigl(\theta,\theta'\bigr)=
\frac{2c^2}{n}\bigl\llvert S\Delta S'\bigr\rrvert \bigl(n-\bigl
\llvert S\Delta S'\bigr\rrvert \bigr),
\]
where $A\Delta B$ is the symmetric difference defined as $(A\cap
B^c)\cup(A^c\cap B)$. By viewing $\llvert  S\Delta S'\rrvert $
as the Hamming
distance of the corresponding indicator functions of the sets, we can
use the Varshamov--Gilbert bound (Lemma~\ref{lemVG}) to pick
$S_1,\ldots,S_N\subset\mathcal{S}$ satisfying
\[
\tfrac{1}{4}n\leq\llvert S_i\Delta S_j\rrvert
\leq\tfrac
{3}{4}n\qquad\mbox{for }i\neq j\in[N],
\]
with $N\geq\exp(c_1n)$, for some $c_1>0$. Hence, we have
\[
\mathcal{M}(\varepsilon,T,\rho)\geq N\geq\exp(c_1n)\qquad\mbox{with }
\varepsilon^2=\frac{c^2}{8n}.
\]
Applying (\ref{eqfanoKL}) of Proposition~\ref{propfano}, we have
\[
\inf_{\hat{\theta}}\sup_{\theta\in\Theta_2}\mathbb{P} \biggl\{
\frac{1}{n^2}\sum_{i,j\in[n]}(\hat{\theta}_{ij}-
\theta _{ij})^2\geq\frac{c^2}{32n} \biggr\}\geq1-
\frac{8c^2n+\log
2}{c_1n}\geq0.8,
\]
where the last inequality holds by choosing a sufficiently small $c$.
Note that the above derivation ignores the fact that $\theta_{ii}=0$
for $i\in[n]$ for the sake of clear presentation. The argument can be
easily made rigorous with slight modification.
Thus, we prove the lower bound for a finite $k$. For $k$ growing with
$n$, a more delicate construction is stated in Section~\ref{secpflower}.

%s3.4 #&#
\subsection{Application to link prediction}

An important application of Theorems~\ref{teomain2} and~\ref
{teomain1} is link prediction.
The link prediction or the network completion problem \cite
{guimera2009missing,lu2011link,zhang2014degree} has practical significances.
Instead of observing the whole adjacency matrix, we observe $\{
A_{ij}:(i,j)\in\Omega\}$ for some $\Omega\subset[n]\times[n]$. The
goal is to infer the unobserved edges.
One example is the biological network. Scientific study showed that
only 80\% of the molecular interactions in cells of Yeast are known
\cite{yu2008high}. Accurate prediction of those unseen interactions
can greatly reduce the costs of biological experiments.
To tackle the problem of link prediction, we consider
a modification of the constrained least square program, which is
defined~as
%
%e3.1 #&#
%
\begin{equation}
\min\llVert \theta\rrVert ^2-\frac{2n^2}{\llvert \Omega
\rrvert }\sum
_{(i,j)\in\Omega
}A_{ij}\theta_{ij}\qquad\mbox{s.t. }
\theta\in\Theta_k.\label
{eqcompletion}
\end{equation}
The estimator $\hat{\theta}$ obtained from solving (\ref
{eqcompletion}) takes advantage of the underlying block structure of
the network, and is an extension to (\ref{eqestimator}). The number
$\hat{\theta}_{ij}$ can be interpreted as how likely there is an edge
between $i$ and $j$. To analyze the theoretical performance of (\ref
{eqcompletion}), let us assume the set $\Omega$ is obtained by
uniformly sampling with replacement from all edges. In other words,
$\Omega$ may contain some repeated elements.

%th3.4 #&#
%
\begin{teo}\label{teolink-pred}
Assume $\llvert \Omega\rrvert /n^2\geq c$ for a constant
$c\in(0,1]$. For any
constant $C'>0$, there exists some constant $C>0$ only depending on
$C'$ and $c$ such that
\[
\frac{1}{n^2}\sum_{i,j\in[n]}(\hat{
\theta}_{ij}-\theta _{ij})^2\leq C \biggl(
\frac{k^2}{n^2}+\frac{\log k}{n} \biggr),
\]
with probability at least $1-\exp(-C'n\log k)$ uniformly over $\theta
\in\Theta_k$ for all $k\in[n]$.
\end{teo}

The result of Theorem~\ref{teolink-pred} assumes $\llvert \Omega
\rrvert /n^2\geq
c$. For example, when $\llvert \Omega\rrvert /n^2=1/2$, we
only observe at most half
of the edges. Theorem~\ref{teolink-pred} gives rate-optimal link
prediction of the rest of the edges.
In contrast, the low-rank matrix completion approach, though
extensively studied and applied in literature, only gives a rate $k/n$,
which is inferior to that of Theorem~\ref{teolink-pred}.

In the case where the assumption of stochastic block model is not
natural \cite{sarkar2012nonparametric}, we may consider a more general
class of networks generated by a smooth graphon.
This is also a useful assumption to do link prediction. Using the same
estimator (\ref{eqcompletion}) with $k=\lceil n^{\afrac{1}{\alpha
\wedge1+1}} \rceil$, we can obtain the error
\[
\frac{1}{n^2}\sum_{i,j\in[n]}(\hat{
\theta}_{ij}-\theta _{ij})^2\leq C
\biggl(n^{-\afrac{2\alpha}{2\alpha+1}}+\frac{\log
n}{n} \biggr),
\]
with probability at least $1-\exp(-C'n)$ uniformly over $f\in\mathcal
{F}_{\alpha}(M)$ and $\mathbb{P}_{\xi}$, which extends Theorem~\ref
{teomain1}. The proof of Theorem~\ref{teolink-pred} is nearly
identical to that of Theorem~\ref{teomain2} and is omitted in the paper.

%s3.5 #&#
\subsection{Minimax rate for operator norm}

The minimax rates in the paper are all studied under the $\ell_2$
norm, which is the Frobenius norm for a matrix. It is also interesting
to investigate the minimax rate under the matrix operator norm. Recall
that for a matrix $U$, its operator norm $\llVert  U\rrVert
_{\mathrm{op}}$ is the largest
singular value.

%th3.5 #&#
%
\begin{teo}\label{teooperator}
For the stochastic block model $\Theta_k$ with $k\geq2$, we have
\[
\inf_{\hat{\theta}}\sup_{\theta\in\Theta_k}\mathbb{E}\llVert \hat{
\theta}-\theta\rrVert _{\mathrm{op}}^2\asymp n.
\]
\end{teo}

Interestingly, the result of Theorem~\ref{teooperator} does not
depend on $k$ as long as $k\geq2$. The optimal estimator is the
adjacency matrix itself $\hat{\theta}=A$, whose bound under the
operator norm can be derived from standard random matrix theory \cite
{vershynin10}. The lower bound is directly implied from Theorem~\ref
{teoSBMlower} by the following argument:
%
%e3.2 #&#
%
\begin{eqnarray}\label{eqoperator}
\nonumber
\inf_{\hat{\theta}}\sup_{\theta\in\Theta_k}\mathbb {E}
\llVert \hat{\theta}-\theta\rrVert _{\mathrm
{op}}^2&\gtrsim&\inf
_{\hat{\theta}}\sup_{\theta\in\Theta_2}\mathbb{E}\llVert \hat{
\theta}-\theta \rrVert _{\mathrm{op}}^2
\nonumber\\[-8pt]\\[-8pt]\nonumber
&\gtrsim& \inf_{\hat{\theta}\in\Theta_2}\sup_{\theta\in\Theta_2}
\mathbb{E}\llVert \hat{\theta}-\theta \rrVert _{\mathrm{op}}^2
\gtrsim\inf_{\hat{\theta}}\sup_{\theta\in\Theta_2}\mathbb {E}\llVert
\hat{\theta}-\theta\rrVert ^2.
\end{eqnarray}
The first inequality is because $\Theta_2$ is a smaller model than
$\Theta_k$ for $k\geq2$. The second inequality is because of the fact
that we can always project the estimator into the parameter space
without compromising the convergence rate. Then, for $\hat{\theta
},\theta\in\Theta_2$, $\hat{\theta}-\theta$ is a matrix with rank
at most $4$, and we have\vspace*{1pt} the inequality $\llVert \hat{\theta
}-\theta
\rrVert ^2\leq4\llVert \hat{\theta}-\theta\rrVert
_{\mathrm{op}}^2$, which gives the last
inequality. Finally, $\inf_{\hat{\theta}}\sup_{\theta\in\Theta
_2}\mathbb{E}\llVert \hat{\theta}-\theta\rrVert ^2\gtrsim
n$ by Theorem~\ref
{teoSBMlower} implies the desired conclusion.

Theorem~\ref{teooperator} suggests that estimating $\theta$ under
the operator norm is not a very interesting problem, because the
estimator does not need to take advantage of the structure of the space
$\Theta_k$. Due to recent advances in community detection, a more suitable
parameter space for the problem is $\Theta(\beta)\cap\Theta_k$, where
\[
\Theta(\beta)= \Bigl\{\theta=\theta^T=\{\theta_{ij}\}\in
[0,1]^{n\times n}: \theta_{ii}=0, \max_{ij}{
\theta_{ij}}\leq\beta \Bigr\}.
\]
The parameter $\beta$ is understood to be the sparsity of the network
because a smaller $\beta$ leads to less edges of the graph.

%th3.6 #&#
%
\begin{teo}\label{teosparse-op}
For $n^{-1}\leq\beta\leq1$ and $k\geq2$, we have
\[
\inf_{\hat{\theta}}\sup_{\theta\in\Theta(\beta)\cap\Theta
_k}\mathbb{E}\llVert \hat{
\theta}-\theta\rrVert _{\mathrm
{op}}^2\asymp\inf
_{\hat
{\theta}}\sup_{\theta\in\Theta(\beta)}\mathbb{E}\llVert \hat {
\theta}-\theta\rrVert _{\mathrm{op}}^2\asymp\beta n.
\]
\end{teo}

The lower bound of Theorem~\ref{teosparse-op} can be obtained in a
similar way by combining the argument in (\ref{eqoperator}) and a
modified version of Theorem~\ref{teoSBMlower} (see the supplementary
material \cite{supp}). When $\beta\geq n^{-1}\log n$, the upper bound
is still achieved by the adjacency matrix, as is proved in Theorem 5.2
of \cite{lei13}. For $n^{-1}\leq\beta< n^{-1}\log n$, one needs to
replace the rows and columns that have high degrees by zeros in~$A$,
and the upper bound is achieved by this trimmed adjacency matrix. This
was recently established in \cite{chin2015stochastic}.

%s3.6 #&#
\subsection{Relation to community detection}

Community detection is another important problem in network analysis.
The parameter estimation result established in this paper has some
consequences in community detection, especially for the results under
the operator norm in Theorems~\ref{teooperator}~and~\ref
{teosparse-op}. Recent works in community detection \cite
{lei13,chin2015stochastic} show that the bound for $\llVert \hat
{\theta}-\theta\rrVert _{\mathrm{op}}^2$ can be used to derive
the misclassification error
of spectral clustering algorithm applied on the matrix $\hat{\theta
}$. Recall that the spectral clustering algorithm applies $k$-means to
the leading singular vectors of the matrix $\hat{\theta}$. Theorem
\ref{teooperator} justifies the use of adjacency matrix as $\hat
{\theta}$ in spectral clustering because of its minimax optimality
under the operator norm.
Moreover, when the network is in a sparse regime with $n^{-1}\leq\beta
<n^{-1}\log n$, \cite{chin2015stochastic} suggests to use the trimmed
adjacency matrix as $\hat{\theta}$ for spectral clustering. According
to Theorem~\ref{teosparse-op}, the trimmed adjacency matrix is an
optimal estimator of $\theta$ under the operator norm.

On the other hand, the connection between the minimax rates under the
$\ell_2$ norm and community detection is not that close. We illustrate
this point by the case when $k=2$. Let us consider $\theta\in\Theta
_2$, then $\theta_{ij}=Q_{z(i)z(j)}$ for some $2\times2$ symmetric
matrix $Q$ and $z$ is the label function. Suppose the within community
connection probability is greater than the between community connection
probability by a margin of $s$. Namely, assume $Q_{11}\wedge
Q_{22}-Q_{12}\geq s>0$. Then,\vspace*{1pt} for the estimator $\hat{\theta
}_{ij}=\hat{Q}_{\hat{z}(i)\hat{z}(j)}$ with error $\frac
{1}{n^2}\sum_{i,j\in[n]}(\hat{\theta}_{ij}-\theta_{ij})^2\leq
\varepsilon^2$, the number of mis-clustered nodes under $\hat{z}$ is
roughly bounded by $O ((n\varepsilon/s)^2 )$. This is because
when two nodes that have the same labels under $z$ are clustered into
different communities or when two nodes belong to different communities
are clustered into the same one, an estimation error of $O(s^2)$ must
occur. Conversely, bounds on community detection can lead to an
improved bound for parameter estimation. Specifically, when $
(\sqrt{Q_{11}\wedge Q_{22}}-\sqrt{Q_{12}} )^2>2n^{-1}\log n$
and $\llvert  z^{-1}(1)\rrvert =\llvert  z^{-1}(2)\rrvert =n/2$, \cite
{mossel2014consistency,hajek2014achieving} show that there exists a
strongly consistent estimator of $z$ in the sense that the
misclassification error is $0$ with high probability. In this case, the
estimation error of $\theta$ under the loss $\frac{1}{n^2}\sum_{i,j\in[n]}(\hat{\theta}_{ij}-\theta_{ij})^2$ can be improved to
$n^{-2}$ from $n^{-1}$.

Generally, parameter estimation and community detection are different
problems of network analysis. When $\{Q_{ab}\}_{a,b\in[k]}$ all take
the same value, it is impossible to do community detection, but
parameter estimation would be easy. Thus, good parameter estimation
result does not necessarily imply consistent community detection.
General minimax rates of the community detection problem are recently
established in \cite{yezhang15,gao15}.

%s4 #&#
\section{Proofs} \label{secproof}

We present the proofs of the main results in this section. The upper
bounds Theorems~\ref{teomain2}~and~\ref{teomain1} are proved
in Section~\ref{secpfupper}. The lower bound Theorem~\ref
{teoSBMlower} is proved in Section~\ref{secpflower}.

%s4.1 #&#
\subsection{Proofs of Theorems \texorpdfstring{\protect\ref{teomain2}}{2.1} and \texorpdfstring{\protect\ref{teomain1}}{2.3}} \label{secpfupper}

This section is devoted to proving the upper bounds. We first prove
Theorem~\ref{teomain2} and then prove Theorem~\ref{teomain1}.

Let us first give an outline of the proof of Theorem~\ref{teomain2}.
In the definition of the class $\Theta_k$, we denote the true value on
each block by $\{Q_{ab}^*\}\in[0,1]^{k\times k}$ and the oracle
assignment by $z^*\in\mathcal{Z}_{n,k}$ such that $\theta
_{ij}=Q^*_{z^*(i)z^*(j)}$ for any $i\neq j$. To facilitate the proof,
we introduce the following notation. For the estimated $\hat{z}$,
define $\{\tilde{Q}_{ab}\}\in[0,1]^{k\times k}$ by $\tilde
{Q}_{ab}=\bar{\theta}_{ab}(\hat{z})$, and also define $\tilde
{\theta}_{ij}=\tilde{Q}_{\hat{z}(i)\hat{z}(j)}$ for any $i\neq j$.
The diagonal elements $\{\tilde{\theta}_{ii}\}$ are defined as zero
for all $i\in[n]$. By the definition of the estimator (\ref
{eqestimator}), we have
\[
L(\hat{Q},\hat{z})\leq L\bigl(Q^*,z^*\bigr),
\]
which can be rewritten as
%
%e4.1 #&#
%
\begin{equation}
\llVert \hat{\theta}-A\rrVert ^2\leq\llVert \theta -A\rrVert
^2. \label{eqbasic}
\end{equation}
The left-hand side of (\ref{eqbasic}) can be decomposed as
%
%e4.2 #&#
%
\begin{equation}
\llVert \hat{\theta}-\theta\rrVert ^2+2 \langle\hat {\theta}-\theta,
\theta -A \rangle+\llVert \theta -A\rrVert ^2. \label{eqleftside}
\end{equation}
Combining (\ref{eqbasic}) and (\ref{eqleftside}), we have
%
%e4.3 #&#
%
\begin{equation}
\llVert \hat{\theta}-\theta\rrVert ^2\leq2 \langle \hat{\theta}-
\theta, A-\theta \rangle. \label{eqrightside}
\end{equation}
The right-hand side of (\ref{eqrightside}) can be bounded as
%
%e4.4 #&#
%e4.5 #&#
%
\begin{eqnarray}
\nonumber
\langle\hat{\theta}-\theta, A-\theta \rangle &=& \langle\hat {\theta}-
\tilde{\theta}, A-\theta \rangle + \langle\tilde{\theta }-\theta, A-\theta \rangle
\\
\label{eqTBC} &\leq& \llVert \hat{\theta}-\tilde{\theta }\rrVert \biggl\llvert
\biggl\langle\frac{\hat{\theta}-\tilde
{\theta}}{\llVert \hat{\theta }-\tilde{\theta}\rrVert }, A-\theta \biggr\rangle\biggr\rrvert
\\
\label{eqTBC1} &&{} + \bigl(\llVert \tilde{\theta}-\hat{\theta }\rrVert +\llVert
\hat {\theta}-\theta\rrVert \bigr)\biggl\llvert \biggl\langle\frac
{\tilde{\theta }-\theta}{\llVert \tilde{\theta}-\theta\rrVert },
A-\theta \biggr\rangle\biggr\rrvert.
\end{eqnarray}
Using Lemmas~\ref{lemaverage}--\ref{lempartition2}, the following
three terms:
%
%e4.6 #&#
%
\begin{equation}
\llVert \hat{\theta}-\tilde{\theta}\rrVert,\qquad\biggl\llvert \biggl\langle
\frac{\hat {\theta}-\tilde{\theta}}{\llVert \hat{\theta}-\tilde{\theta }\rrVert }, A-\theta \biggr\rangle\biggr\rrvert,\qquad\biggl\llvert \biggl
\langle\frac{\tilde
{\theta }-\theta}{\llVert \tilde{\theta}-\theta\rrVert }, A-\theta \biggr\rangle\biggr\rrvert \label{eq3terms}
\end{equation}
can all be bounded by $C\sqrt{k^2+n\log k}$ with probability at least
\[
1-3\exp\bigl(-C'n\log k\bigr).
\]
Combining these bounds with (\ref{eqTBC}),
(\ref{eqTBC1}) and (\ref{eqrightside}), we get
\[
\llVert \hat{\theta}-\theta\rrVert ^2\leq C_1
\bigl(k^2+k\log n \bigr),
\]
with probability at least $1-3\exp(-C'n\log k)$. This gives the
conclusion of Theorem~\ref{teomain2}. The details of the proof is
stated in the later part of the section. To prove Theorem~\ref
{teomain1}, we use Lemma~\ref{lembias} to approximate the
nonparametric graphon by the stochastic block model. With similar
arguments above, we get
\[
\llVert \hat{\theta}-\theta\rrVert ^2\leq C_2
\bigl(k^2+k\log n+ n^2 k^{-2(\alpha\wedge1)} \bigr),
\]
with high probability. Choosing the best $k$ gives the conclusion of
Theorem~\ref{teomain1}.

Before stating the complete proofs, let us first present the following
lemmas, which bound the three terms in (\ref{eq3terms}),
respectively. The proofs of the lemmas will be given in the
supplementary material \cite{supp}.

%le4.1 #&#
%
\begin{lemma} \label{lemaverage}
For any constant $C'>0$, there exists a constant $C>0$ only depending
on $C'$, such that
\[
\llVert \hat{\theta}-\tilde{\theta}\rrVert \leq C\sqrt {k^2+n\log
k},
\]
with probability at least $1-\exp(-C'n\log k)$.
\end{lemma}

%le4.2 #&#
%
\begin{lemma} \label{lempartition1}
For any constant $C'>0$, there exists a constant $C>0$ only depending
on $C'$, such that
\[
\biggl\llvert \biggl\langle\frac{\tilde{\theta}-\theta}{\llVert
\tilde{\theta }-\theta\rrVert }, A-\theta \biggr\rangle \biggr
\rrvert \leq C\sqrt{n\log k},
\]
with probability at least $1-\exp(-C'n\log k)$.
\end{lemma}

%le4.3 #&#
%
\begin{lemma} \label{lempartition2}
For any constant $C'>0$, there exists a constant $C>0$ only depending
on $C'$, such that
\[
\biggl\llvert \biggl\langle\frac{\hat{\theta}-\tilde{\theta}}{\llVert \hat{\theta }-\tilde{\theta}\rrVert }, A-\theta \biggr\rangle\biggr
\rrvert \leq C\sqrt{k^2+n\log k},
\]
with probability at least $1-\exp(-C'n\log k)$.
\end{lemma}

\begin{pf*}{Proof of Theorem~\ref{teomain2}}
Combining the bounds for (\ref{eq3terms}) with (\ref{eqTBC}), (\ref
{eqTBC1}) and (\ref{eqrightside}), we have
\[
\llVert \hat{\theta}-\theta\rrVert ^2\leq2C\llVert \hat{\theta}-
\theta\rrVert \sqrt {k^2+n\log k}+4C^2
\bigl(k^2+n\log k \bigr),
\]
with probability at least $1-3\exp(-C'n\log k)$. Solving the above
equation, we get
\[
\llVert \hat{\theta}-\theta\rrVert ^2\leq C_1
\bigl(k^2+n\log k \bigr),
\]
with probability at least $1-3\exp(-C'n\log k)$. This proves the high
probability bound. To get the bound in expectation, we use the
following inequality:
\begin{eqnarray*}
&& \mathbb{E}n^{-2}\llVert \hat{\theta}-\theta\rrVert ^2
\\
&&\qquad \leq \mathbb{E} \bigl(n^{-2}\llVert \hat{\theta}-\theta \rrVert
^2\mathbb {I}\bigl\{n^{-2}\llVert \hat{\theta}-\theta
\rrVert ^2\leq \varepsilon^2\bigr\} \bigr)
\\
&&\quad\qquad{} + \mathbb{E}
\bigl(n^{-2}\llVert \hat{\theta}-\theta\rrVert ^2\mathbb{I}
\bigl\{ n^{-2}\llVert \hat{\theta}-\theta\rrVert ^2>
\varepsilon^2\bigr\} \bigr)
\\
&&\qquad \leq \varepsilon^2 + \mathbb{P} \bigl(n^{-2}\llVert \hat{
\theta }-\theta \rrVert ^2>\varepsilon^2 \bigr) \leq
\varepsilon^2+3\exp\bigl(-C'n\log k\bigr),
\end{eqnarray*}
where $\varepsilon^2=C_1 (\frac{k^2}{n^2}+\frac{\log k}{n}
)$. Since $\varepsilon^2$ is the dominating term, the proof is complete.
\end{pf*}

To prove Theorem~\ref{teomain1}, we need to redefine $z^*$ and $Q^*$.
We choose $z^*$ to be the one used in Lemma~\ref{lembias}, which
implies a good approximation of $\{\theta_{ij}\}$ by the stochastic
block model. With this $z^*$, define $Q^*$ by letting $Q_{ab}^*=\bar
{\theta}_{ab}(z^*)$ for any $a,b\in[k]$. Finally, we define $\theta
_{ij}^*=Q^*_{z^*(i)z^*(j)}$ for all $i\neq j$. The diagonal elements
$\theta_{ii}^*$ are set as zero for all $i\in[n]$. Note that for the
stochastic block model, we have $\theta=\theta^*$. The proof of
Theorem~\ref{teomain1} requires another lemma.
%le4.4 #&#

\begin{lemma} \label{lempartition3}
For any constant $C'>0$, there exists a constant $C>0$ only depending
on $C'$, such that
\[
\biggl\llvert \biggl\langle\frac{\tilde{\theta}-\theta^*}{\llVert \tilde{\theta }-\theta^*\rrVert }, A-\theta \biggr\rangle\biggr
\rrvert \leq C\sqrt{n\log k},
\]
with probability at least $1-\exp(-C'n\log k)$.
\end{lemma}

The proof of Lemma~\ref{lempartition3} is identical to the proof of
Lemma~\ref{lempartition1}, and will be omitted in the paper.

\begin{pf*}{Proof of Theorem~\ref{teomain1}}
Using the similar argument as outlined in the beginning of this
section, we get
\[
\bigl\llVert \hat{\theta}-\theta^*\bigr\rrVert ^2\leq2 \bigl\langle
\hat{\theta}-\theta ^*, A-\theta^* \bigr\rangle,
\]
whose right-hand side can be bounded as
\begin{eqnarray*}
&& \bigl\langle\hat{\theta}-\theta^*, A-\theta^* \bigr\rangle
\\
&&\qquad = \langle\hat{\theta}-\tilde{\theta}, A-\theta \rangle+ \bigl\langle\tilde {
\theta}-\theta^*, A-\theta \bigr\rangle+ \bigl\langle\hat{\theta}-\theta ^*,
\theta-\theta^* \bigr\rangle
\\
&&\qquad \leq \llVert \hat{\theta}-\tilde{\theta}\rrVert \biggl\llvert \biggl\langle
\frac{\hat {\theta}-\tilde{\theta}}{\llVert \hat{\theta}-\tilde{\theta }\rrVert }, A-\theta \biggr\rangle\biggr\rrvert + \bigl(\llVert \tilde{
\theta}-\hat {\theta }\rrVert +\bigl\llVert \hat{\theta}-\theta^*\bigr\rrVert
\bigr)\biggl\llvert \biggl\langle\frac{\tilde {\theta}-\theta^*}{\llVert \tilde{\theta}-\theta^*\rrVert }, A-\theta \biggr\rangle\biggr
\rrvert
\\
&&\quad\qquad{}  + \bigl\llVert \hat{\theta}-\theta^*\bigr\rrVert \bigl\llVert \theta-\theta^*
\bigr\rrVert.
\end{eqnarray*}
To better organize what we have obtained, let us introduce the notation
\[
L = \bigl\llVert \hat{\theta}-\theta^*\bigr\rrVert,\qquad R = \llVert \tilde{
\theta}-\hat{\theta}\rrVert, \qquad B = \bigl\llVert \theta-\theta^*\bigr\rrVert,
\]
\[
E = \biggl\llvert \biggl\langle\frac{\hat{\theta}-\tilde{\theta
}}{\llVert \hat {\theta}-\tilde{\theta}\rrVert }, A-\theta \biggr\rangle\biggr
\rrvert, \qquad F = \biggl\llvert \biggl\langle\frac{\tilde{\theta}-\theta^*}{\llVert \tilde{\theta }-\theta^*\rrVert }, A-\theta \biggr
\rangle\biggr\rrvert. %
\]
Then, by the derived inequalities, we have
\[
L^2\leq2RE+2(L+R)F+2LB.
\]
It can be rearranged as
\[
L^2\leq2(F+B)L+2(E+F)R.
\]
By solving this quadratic inequality of $L$, we can get
%
%e4.7 #&#
%
\begin{equation}
L^2\leq\max\bigl\{ 16(F+B)^2, 4R(E+F)\bigr\}.
\label{eqboundL}
\end{equation}
By Lemma~\ref{lembias}, Lemma~\ref{lemaverage}, Lemma~\ref
{lempartition2} and Lemma~\ref{lempartition3}, for any constant
$C'>0$, there exist constants $C$ only depending on $C', M$, such that
\begin{eqnarray*}
B^2 &\leq& Cn^2 \biggl(\frac{1}{k^2}
\biggr)^{\alpha\wedge1}, \qquad F^2 \leq Cn\log k,
\\
R^2 &\leq& C\bigl(k^2+n\log k\bigr), \qquad E^2
\leq C\bigl(k^2+n\log k\bigr), %
\end{eqnarray*}
with probability at least $1-\exp(-C'n)$. By (\ref{eqboundL}), we have
%
%e4.8 #&#
%
\begin{equation}
\label{equL} L^2\leq C_1 \biggl(n^2
\biggl(\frac{1}{k^2} \biggr)^{\alpha\wedge
1}+k^2+n\log k \biggr)
\end{equation}
with probability at least $1-\exp(-C'n)$ for some constant $C_1$.
Hence, there is some constant $C_2$ such that
\begin{eqnarray*}
\frac{1}{n^2}\sum_{ij}(\hat{
\theta}_{ij}-\theta_{ij})^2 &\leq&
\frac{2}{n^2} \bigl(L^2+B^2 \bigr)
\\
&\leq& C_2 \biggl( \biggl(\frac{1}{k^2} \biggr)^{\alpha\wedge1}+
\frac
{k^2}{n^2}+\frac{\log k}{n} \biggr),
\end{eqnarray*}
with probability at least $1-\exp(-C'n)$. When $\alpha\geq1$, we
choose $k=\lceil\sqrt{n} \rceil$, and the bound is
$C_3n^{-1}\log n$
for some constant $C_3$ only depending on $C'$ and $M$. When $\alpha
<1$, we choose $k=\lceil n^{\afrac{1}{\alpha+1}} \rceil$. Then the
bound is
$C_4n^{-\afrac{2\alpha}{\alpha+1}}$
for some constant $C_4$ only depending on $C'$ and $M$. This completes
the proof.
\end{pf*}

%s4.2 #&#
\subsection{Proof of Theorem \texorpdfstring{\protect\ref{teoSBMlower}}{2.2}} \label{secpflower}

This section is devoted to proving the lower bounds.
For any probability measures $\mathbb{P},\mathbb{Q}$, define the
Kullback--Leibler divergence by $D(\mathbb{P}\|
\mathbb{Q})=\int
(\log\frac{d\mathbb{P}}{d\mathbb{Q}} )\,d\mathbb{P}$. The
chi-squared divergence is defined by $\chi^2(\mathbb{P}\|
\mathbb
{Q})=\int (\frac{d\mathbb{P}}{d\mathbb{Q}} )\,d\mathbb
{P}-1$. To prove minimax lower bounds, we need the following proposition.

%pr4.1 #&#
%
\begin{proposition} \label{propfano}
Let $(\Theta,\rho)$ be a metric space and $\{\mathbb{P}_{\theta
}:\theta\in\Theta\}$ be a collection of probability measures. For
any totally bounded $T\subset\Theta$, define the Kullback--Leibler
diameter and the chi-squared diameter of $T$ by
\[
d_{\mathrm{KL}}(T)=\sup_{\theta,\theta'\in T}D\bigl(\mathbb{P}_{\theta
}
\llVert \mathbb{P}_{\theta'}\bigr),\qquad d_{\chi^2}(T)=\sup
_{\theta,\theta
'\in T}\chi^2\bigl(\mathbb{P}_{\theta}
\rrVert \mathbb{P}_{\theta'}\bigr).
\]
Then
%
%e4.9 #&#
%e4.10 #&#
%
\begin{eqnarray}
\label{eqfanoKL}
\inf_{\hat{\theta}}\sup_{\theta\in\Theta
}\mathbb{P}_{\theta} \biggl\{\rho^2 \bigl(\hat{\theta}(X),\theta
\bigr)\geq\frac{\varepsilon^2}{4} \biggr\} &\geq & 1-\frac{d_{\mathrm
{KL}}(T)+\log2}{\log\mathcal{M}(\varepsilon,T,\rho)},
\\
\label{eqfanochi2}\qquad \inf_{\hat{\theta}}\sup_{\theta\in\Theta
}\mathbb{P}_{\theta} \biggl\{\rho^2 \bigl(\hat{\theta}(X),\theta
\bigr)\geq\frac{\varepsilon^2}{4} \biggr\} &\geq& 1-\frac{1}{\mathcal
{M}(\varepsilon,T,\rho)}-\sqrt{
\frac{d_{\chi^2}(T)}{\mathcal
{M}(\varepsilon,T,\rho)}},
\end{eqnarray}
for any $\varepsilon>0$.
\end{proposition}

Inequality (\ref{eqfanoKL}) is the classical Fano's inequality.
The version we present here is by \cite{yu97}. Inequality (\ref
{eqfanochi2}) is a generalization of the classical Fano's inequality
by using chi-squared divergence instead of KL divergence. It is due to
\cite{guntuboyina11}. We use it here as an alternative of Assouad's
lemma to get the corresponding in-probability lower bound. In this
section, the parameter is a matrix $\{\theta_{ij}\}\in[0,1]^{n\times
n}$. The metric we consider is
\[
\rho^2\bigl(\theta,\theta'\bigr)=\frac{1}{n^2}\sum
_{ij}\bigl(\theta_{ij}-\theta
_{ij}'\bigr)^2.
\]
Let us give bounds for KL divergence and chi-squared divergence under
random graph model. Let $P_{\theta_{ij}}$ denote the probability of
$\operatorname{Bernoulli}(\theta_{ij})$. Given $\theta=\{\theta
_{ij}\}\in[0,1]^{n\times n}$, the probability $\mathbb{P}_{\theta}$
stands for the product measure $\bigotimes_{i,j\in[n]}P_{\theta_{ij}}$
throughout this section.
%pr4.2 #&#

\begin{proposition} \label{propprobdistance}
For any $\theta, \theta' \in[1/2,3/4]^{n\times n}$, we have
%
%e4.11 #&#
%
\begin{eqnarray}\label{eqprobdistance}
D\bigl(\mathbb{P}_\theta\llVert
\mathbb{P}_{\theta'}\bigr) &\le& 8 \sum_{ij}
\bigl(\theta_{ij} - \theta'_{ij}
\bigr)^2,
\nonumber\\[-8pt]\\[-8pt]\nonumber
\chi^2\bigl(\mathbb{P}_\theta
\rrVert \mathbb{P}_\theta'\bigr) &\le& \exp \biggl(8
\sum_{ij} \bigl(\theta _{ij} -
\theta'_{ij}\bigr)^2 \biggr).
\end{eqnarray}
\end{proposition}

The proposition will be proved in the supplementary material \cite
{supp}. We also need the following Varshamov--Gilbert bound. The
version we present here is due to \cite{massart2007}, Lemma 4.7.
%le4.5 #&#

\begin{lemma}\label{lemVG}
There exists a subset $\{\omega_1,\ldots,\omega_N\}\subset\{0,1\}
^d$ such that
%
%e4.12 #&#
%
\begin{equation}
\rho_{H}(\omega_i, \omega_j) \triangleq
\llVert \omega_i-\omega _j\rrVert ^2\geq
\frac{d}{4}\qquad\mbox{for any }i\neq j\in [N],\label{eqdefH}
\end{equation}
for some $N\geq\exp (d/8)$.
\end{lemma}

\begin{pf*}{Proof of Theorem~\ref{teoSBMlower}}
By the definition of the parameter space $\Theta_k$, we rewrite the
minimax rate as
\begin{eqnarray*}
&& \inf_{\hat{\theta}}\sup_{\theta\in\Theta_k}\mathbb{P} \biggl\{
\frac{1}{n^2}\sum_{ij}(\hat{\theta}_{ij}-
\theta_{ij})^2\geq \varepsilon^2 \biggr\}
\\
&&\qquad = \inf_{\hat{\theta}}\sup_{Q=Q^T \in[0,1]^{k \times k}} \sup
_{z
\in\mathcal{Z}_{n,k}}\mathbb{P} \biggl\{\frac{1}{n^2}\sum
_{i\neq
j}(\hat{\theta}_{ij}-Q_{z(i)z(j)})^2
\geq\varepsilon^2 \biggr\}. \label
{eqtheta}
\end{eqnarray*}
If we fix a $z \in\mathcal{Z}_{n,k}$, it will be direct to derive the
lower bound $k^2/n^2$ for estimating~$Q$. On the other hand, if we fix
$Q$ and let $z$ vary, it will become a new type of convergence rate due
to the unknown label and we name it as the clustering rate, which is at
the order of $n^{-1}\log k$. In the following arguments, we will prove
the two different rates separately and then combine them together to
get the desired in-probability lower bound.

Without loss of generality, we consider the case where both $n/k$ and
$k/2$ are integers. If they are not, let $k'=2{ \lfloor{k/2}
 \rfloor}$ and $n' =
\lfloor n/k' \rfloor k'$. By restricting the unknown parameters to the
smaller class $Q'=(Q')^T \in[0,1]^{k' \times k'}$ and $z' \in\mathcal
{Z}_{n',k'}$, the following lower bound argument works for this smaller
class. Then it also provides a lower bound for the original larger
class.

\subsubsection*{Nonparametric rate} First we fix a $z \in\mathcal
{Z}_{n,k}$. For each $a\in[k]$, we define $z^{-1}(a)= \{
(a-1)n/k+1,\ldots,an/k  \}$.
Let $\Omega=\{0,1\}^d$ be the set of all binary sequences of length
$d=k(k-1)/2$. For any $\omega=\{\omega_{ab}\}_{1\leq b< a\leq k} \in
\Omega$, define a $k \times k$ matrix $Q^{\omega}=(Q_{ab}^{\omega
})_{k \times k}$ by
%
%e4.13 #&#
%
\begin{eqnarray}\label{eqyulu}
Q_{ab}^{\omega}&=&Q_{ba}^{\omega}=
\frac{1}{2} + \frac{c_1k}{n}\omega _{ab}\qquad\mbox{for } a > b \in[k]\quad\mbox{and}
\nonumber\\[-8pt]\\[-8pt]\nonumber
Q_{aa}^{\omega} &=& \tfrac{1}{2}\qquad\mbox{for }a\in[k],
\end{eqnarray}
where $c_1$ is a constant that we are going to specify later. Define
$\theta^{\omega} = (\theta_{ij}^{\omega})_{n \times n}$ with
$\theta_{ij}^{\omega}=Q_{z(i)z(j)}^{\omega}$ for $i\neq j$ and
$\theta_{ii}^{\omega}=0$. The subspace we consider is $T_1=\{\theta
^{\omega}: \omega\in\Omega\} \subset\Theta_k$. To apply (\ref
{eqfanochi2}), we need to upper bound $\sup_{\theta,\theta'\in
T_1}\chi^2(\mathbb{P}_{\theta}\|\mathbb
{P}_{\theta'})$ and lower
bound $\mathcal{M}(\varepsilon,T_1,\rho)$. For any $\theta^{\omega},
\theta^{\omega'} \in T_1$, from (\ref{eqprobdistance}) and (\ref
{eqyulu}), we get
%
%e4.14 #&#
%
\begin{eqnarray}\label{eqrho2H}
\nonumber
\chi^2(\mathbb{P}_{\theta^{\omega}}\|
\mathbb{P}_{\theta
^{\omega'}}) &=& \exp \biggl(8\sum_{i,j\in[n]}
\bigl(\theta_{ij}^{\omega
}-\theta_{ij}^{\omega'}
\bigr)^2 \biggr)
\nonumber\\[-8pt]\\[-8pt]\nonumber
&\le& \exp \biggl( \frac{8n^2}{k^2}\sum
_{a,b\in
[k]}\bigl(Q_{ab}^{\omega} -
Q_{ab}^{\omega'}\bigr)^2 \biggr) \le\exp\bigl(8
c_1^2 k^2\bigr),
\end{eqnarray}
where we choose sufficiently small $c_1$ so that $\theta_{ij}^{\omega
},\theta_{ij}^{\omega'}\in[1/2,3/4]$ is satisfied.
To lower bound the packing number, we reduce the metric $\rho(\theta
^{\omega},\theta^{\omega'})$ to $\rho_H(\omega,\omega')$ defined
in (\ref{eqdefH}). In view of (\ref{eqyulu}), we get
%
%e4.15 #&#
%
\begin{equation}
\label{eqdistlowerbound} \rho^2\bigl(\theta^{\omega},
\theta^{\omega'}\bigr) \ge\frac{1}{k^2} \sum
_{1 \le b < a \le k} \bigl(Q_{ab}^{\omega
}-Q^{\omega'}_{ab}
\bigr)^2 = \frac{c_1^2}{n^2} \rho_{H} \bigl(\omega,
\omega'\bigr).
\end{equation}
By Lemma~\ref{lemVG}, we can find a subset $S \subset\Omega$ that
satisfies the following properties: (a) $\llvert  S\rrvert  \ge
\exp (d/8)$ and (b)
$\rho_{H} (\omega, \omega') \ge d/4$ for any $\omega, \omega' \in
S$. From (\ref{eqdistlowerbound}), we have
\[
\mathcal{M}(\varepsilon,T_1,\rho)\geq\llvert S\rrvert \ge \exp (d/8) =
\exp \bigl(k(k-1)/16\bigr),
\]
with $\varepsilon^2=\frac{c_1k(k-1)}{8n^2}$. By choosing sufficiently
small $c_1$, together with (\ref{eqrho2H}), we get
%
%e4.16 #&#
%
\begin{equation}
\inf_{\hat{\theta}}\sup_{\theta\in T_1}\mathbb{P} \biggl\{
\frac
{1}{n^2}\sum_{ij}(\hat{
\theta}_{ij}-\theta_{ij})^2\geq\frac
{C_1k^2}{n^2}
\biggr\}\geq0.9, \label{eqlowerS1}
\end{equation}
by (\ref{eqfanochi2}) for sufficiently large $k$ with some constant
$C_1>0$. When $k$ is not sufficiently large, that is, $k\leq O(1)$,
then it is easy to see that $n^{-2}$ is always the correct order of
lower bound. Since $n^{-2}\asymp k^2/n^2$ when $k\leq O(1)$, $k^2/n^2$
is also a valid lower bound for small $k$.

\subsubsection*{Clustering rate} We are going to fix a $Q$ that has
the following form:
%
%e4.17 #&#
%
\begin{equation}
Q=\lleft[\matrix{ 0 & B
\cr
B^T & 0} \rright],
\label{eqspecialQ}
\end{equation}
where $B$ is a $(k/2)\times(k/2)$ matrix. By Lemma~\ref{lemVG}, when
$k$ is\vspace*{1pt} sufficiently large, we can find $\{\omega_1,\ldots, \omega
_{k/2}\} \subset\{0,1\}^{k/2}$ such that $\rho_H(\omega_a,\omega_b)
\ge k/8$ for all $a\neq b \in[k/2]$. Fixing such $\{\omega_1,\ldots,\omega_{k/2}\}$, define $B=(B_1,B_2,\ldots,B_{k/2})$ by letting
$B_a=\frac{1}{2} + \sqrt{\frac{c_2 \log k}{n}}\omega_a$ for $a\in
[k/2]$. With\vspace*{1pt} such construction, it is easy to see that for any $a\neq
b\in[k/2]$,
%
%e4.18 #&#
%
\begin{equation}
\llVert B_a-B_b\rrVert ^2\geq
\frac{c_2k\log k}{8n}. \label
{eqBspecial}
\end{equation}
Define a subset of $\mathcal{Z}_{n,k}$ by
\begin{eqnarray*}
\mathcal{Z}&=& \biggl\{ z \in\mathcal{Z}_{n,k}: \bigl\llvert
z^{-1}(a)\bigr\rrvert =\frac
{n}{k}\mbox{ for } a \in[k],
\\
&&{}
z^{-1}(a)= \biggl\{\frac{(a-1)n}{k}+1,\ldots,\frac
{an}{k}
\biggr\}\mbox{ for } a \in[k/2] \biggr\}.
\end{eqnarray*}
For each $z\in\mathcal{Z}$, define $\theta^z$ by $\theta
_{ij}^z=Q_{z(i)z(j)}$ for $i\neq j$ and $\theta_{ii}^z=0$.
The subspace we consider is $T_2=\{\theta^z: z\in\mathcal{Z}\}
\subset\Theta_{n,k}$. To apply (\ref{eqfanoKL}), we need to upper
bound $\sup_{\theta,\theta\in T_2}D(\mathbb{P}_\theta\|
\mathbb
{P}_{\theta'})$ and lower bound $\log\mathcal{M}(\varepsilon,T_2,\rho
)$. By (\ref{eqprobdistance}), for any \mbox{$\theta,\theta' \in T_2$},
%
%e4.19 #&#
%
\begin{equation}
\label{equKLd} D(\mathbb{P}_\theta\|
\mathbb{P}_{\theta'}) \le8 \sum_{ij} \bigl(
\theta_{ij}-\theta'_{ij}\bigr)^2
\le8n^2 c_2 \frac{\log k}{n} = 8c_2 n \log k.
\end{equation}
Now we are going to give a lower bound of the packing number $\log
\mathcal{M}(\varepsilon,T_2,\rho)$ with $\varepsilon^2=(c_2\log k)/(48n)$
for the $c_2$ in (\ref{eqBspecial}). Due to the construction of $B$,
there is a one-to-one correspondence between $T_2$ and $\mathcal{Z}$.
Thus, $\log\mathcal{M}(\varepsilon,T_2,\break \rho)=\log\mathcal{M}
(\varepsilon, \mathcal{Z},\rho_1)$ for some metric $\rho_1$ on
$\mathcal{Z}$ defined by $\rho_1(z,w)=\rho(\theta^z, \theta^w)$.
Given any $z\in\mathcal{Z}$, define its $\varepsilon$-neighborhood by
$B(z, \varepsilon) = \{ w \in\mathcal{Z}: \rho_1(z,w) \le\varepsilon\}
$. Let $S$ be the packing set in $\mathcal{Z}$ with cardinality
$\mathcal{M}(\varepsilon,\mathcal{Z},\rho_1)$. We claim that $S$ is
also the covering set of $\mathcal{Z}$ with radius $\varepsilon$,
because otherwise there is some point in $\mathcal{Z}$ which is at
least $\varepsilon$ away from every point in $S$, contradicting the
definition of $\mathcal{M}(\varepsilon,\mathcal{Z},\rho_1)$.
This implies the fact $\bigcup_{z \in S} B (z, \varepsilon) = \mathcal
{Z}$, which leads to
\[
\llvert \mathcal{Z}\rrvert \le\sum_{z \in S} \bigl
\llvert B (z, \varepsilon)\bigr\rrvert \le\llvert S\rrvert \max_{z
\in S}
\bigl\llvert B (z, \varepsilon)\bigr\rrvert.
\]
Thus, we have
%
%e4.20 #&#
%
\begin{equation}
\mathcal{M}(\varepsilon,\mathcal{Z},\rho_1)=\llvert S\rrvert \geq
\frac{\llvert \mathcal
{Z}\rrvert }{\max_{z \in S} \llvert  B (z, \varepsilon)\rrvert }. \label{eqVRatio}
\end{equation}
Let us upper bound $\max_{z \in S} \llvert  B (z, \varepsilon)\rrvert $ first.
For any $z,w \in\mathcal{Z}$, by the construction of $\mathcal{Z}$,
$z(i)=w(i)$ when $i\in[n/2]$ and $\llvert  z^{-1}(a)\rrvert
=n/k$ for each $a\in
[k]$. Hence,
\begin{eqnarray*}
\rho_1^2(z,w) &\geq& \frac{1}{n^2} \sum
_{1 \le i \le n/2 <j \le n} (Q_{z(i)z(j)} - Q_{w(i)w(j)})^2
\\
&=& \frac{1}{n^2} \sum_{n/2 < j \le n} \sum
_{1 \le a \le k/2} \sum_{i \in z^{-1}(a)}
(Q_{az(j)} - Q_{aw(j)})^2
\\
&=& \frac{1}{n^2} \sum_{n/2< j \le n} \frac{n}{k}
\llVert B_{z(j)} - B_{w(j)}\rrVert ^2
\\
&\ge& \frac{c_2 \log k}{8n^2} \bigl\llvert \bigl\{j: w(j) \neq z(j) \bigr\} \bigr
\rrvert,
\end{eqnarray*}
where the last inequality is due to (\ref{eqBspecial}).
Then for any $w \in B(z, \varepsilon)$, $\llvert \{j: w(j) \neq z(j)\}
\rrvert  \le n/6$
under the choice $\varepsilon^2=(c_2\log k)/(48n)$. This implies
\[
\bigl\llvert B(z,\varepsilon)\bigr\rrvert \leq{n\choose n/6}k^{n/6}\leq
(6e)^{n/6}k^{n/6}\leq \exp \biggl(\frac{1}{4}n\log k
\biggr).
\]
Now we lower bound $\llvert \mathcal{Z}\rrvert $.
Note that by Stirling's formula
\[
\llvert \mathcal{Z}\rrvert = \frac{(n/2)!}{[(n/k)!]^{k/2}} = \exp \biggl(
\frac
{1}{2 }n\log k + o(n\log k) \biggr) \ge\exp \biggl(\frac{1}{3} n
\log k \biggr).
\]
By (\ref{eqVRatio}), we get $\log\mathcal{M}(\varepsilon,T,\rho
)=\log\mathcal{M}(\varepsilon,\mathcal{Z},\rho_1) \ge(1/12) n\log
k$. Together with (\ref{equKLd}) and using (\ref{eqfanoKL}), we have
%
%e4.21 #&#
%
\begin{equation}
\inf_{\hat{\theta}}\sup_{\theta\in T_2}\mathbb{P} \biggl\{
\frac
{1}{n^2}\sum_{ij}(\hat{
\theta}_{ij}-\theta_{ij})^2\geq\frac
{C_2\log k}{n}
\biggr\}\geq0.9, \label{eqlowerS2}
\end{equation}
with some constant $C_2>0$
for sufficiently small $c_2$ and sufficiently large $k$. When $k$ is
not sufficiently large but $2\leq k\leq O(1)$, the argument in
Section~\ref{secflower} gives the desired lower bound at the order of
$n^{-1}\asymp n^{-1}\log k$. When $k=1$, $n^{-1}\log k=0$ is still a
valid lower bound.

\subsubsection*{Combining the bounds}
Finally, let us
combine (\ref{eqlowerS1}) and (\ref{eqlowerS2}) to get the desired
in-probability lower bound in Theorem~\ref{teoSBMlower} with $C=
(C_1\wedge C_2 )/2$. For any $\theta\in\Theta_k$, by union
bound, we have
\begin{eqnarray*}
&& \mathbb{P} \biggl\{\frac{1}{n^2}\sum_{ij}(
\hat{\theta }_{ij}-\theta_{ij})^2\geq C \biggl(
\frac{k^2}{n^2}+\frac{\log
k}{n} \biggr) \biggr\}
\\
&&\qquad \geq 1-\mathbb{P} \biggl\{\frac{1}{n^2}\sum_{ij}(
\hat{\theta }_{ij}-\theta_{ij})^2\leq
\frac{C_1k^2}{n^2} \biggr\}-\mathbb {P} \biggl\{\frac{1}{n^2}\sum
_{ij}(\hat{\theta}_{ij}-\theta _{ij})^2
\leq\frac{C_2\log k}{n} \biggr\}
\\
&&\qquad = \mathbb{P} \biggl\{\frac{1}{n^2}\sum_{ij}(
\hat{\theta }_{ij}-\theta_{ij})^2\geq
\frac{C_1k^2}{n^2} \biggr\}+\mathbb {P} \biggl\{\frac{1}{n^2}\sum
_{ij}(\hat{\theta}_{ij}-\theta _{ij})^2
\geq\frac{C_2\log k}{n} \biggr\}-1.
\end{eqnarray*}
Taking $\sup$ on both sides, and using the fact $\sup_{z,Q}
(f(z)+g(Q) )=\sup_zf(z)+\sup_Qg(Q)$, we have
\begin{eqnarray*}
&& \sup_{\theta\in\Theta_k}\mathbb{P} \biggl\{\frac{1}{n^2}\sum
_{ij}(\hat{\theta}_{ij}-\theta_{ij})^2
\geq C \biggl(\frac
{k^2}{n^2}+\frac{\log k}{n} \biggr) \biggr\}
\\
&&\qquad \geq \sup_{\theta\in T_1}\mathbb{P} \biggl\{\frac{1}{n^2}\sum
_{ij}(\hat{\theta}_{ij}-
\theta_{ij})^2\geq\frac{C_1k^2}{n^2} \biggr\}
\\
&&\quad\qquad{}+\sup
_{\theta\in T_2}\mathbb{P} \biggl\{\frac{1}{n^2}\sum
_{ij}(\hat{\theta}_{ij}-\theta_{ij})^2
\geq\frac{C_2\log
k}{n} \biggr\}-1,
\end{eqnarray*}
for any estimator $\hat{\theta}$. Plugging the lower bounds (\ref
{eqlowerS1}) and (\ref{eqlowerS2}), we obtain the desired result.
A Markov's inequality argument leads to the lower bound in expectation.
\end{pf*}

%\begin{appendix}
%\section{}
%\end{appendix}

\section*{Acknowledgements}
We want to thank Zongming Ma for helpful discussion on
the relation between graphon estimation and link prediction, and to thank
the Associate Editor and the referee for their constructive comments and
suggestions that lead to the improvement of the paper.

\begin{supplement}[id=suppA]
\stitle{Supplement to ``Rate-optimal graphon estimation.''\\}
\slink[doi]{10.1214/15-AOS1354SUPP} %[doi,text={...}] - jei reikia
%suskaldyti doi
\sdatatype{.pdf}
\sfilename{aos1354\_supp.pdf}
\sdescription{In the supplement, we prove Theorem~\ref
{teographonlower}, Lemmas~\ref{lembias},~\ref{lemaverage},
\ref{lempartition1},~\ref{lempartition2}, Proposition
\ref{propprobdistance} and Theorem~\ref{teosparse-op}.}
\end{supplement}

% imsref loaded by linak, 2015-07-09 16:06:00
% imsref loaded by linak, 2015-07-10 15:31:45

\printaddresses
\end{document}